\documentclass[12pt, reqno]{gen-p-l}
\usepackage{ifthen,amsfonts,amsmath,amssymb,epic,eepic,epsfig,color}

\begingroup

\newtheorem{theorem}{Theorem}[section]
\newtheorem{lemma}[theorem]{Lemma}
\newtheorem{propos}[theorem]{Proposition}

\endgroup

\newtheorem{definition}[theorem]{Definition}

\newtheorem{remark}[theorem]{Remark}

\newtheorem{Con}[theorem]{Conjecture}

\newtheorem{Exa}[theorem]{Example}

\newcommand{\an}{{\rm An}}
\newcommand{\eps}{{\varepsilon}}

\newcommand{\Vol}{{\text{Vol}}}

\newcommand{\diam}{{\text {diam}}}
\newcommand{\dist}{{\text {dist}}}

\def\RR{{\bf  R}}

\def\SS{{\bf  S}}

\def\CC{{\bf C }}

\newcommand{\eqr}[1]{(\ref{#1})}

\newcommand\diff{{\rm Diff}_0}

\newcommand\D{{\mathcal{D}}}

\newcommand\B{{\mathcal{B}}}
\newcommand\va{{\mathcal{V} (M)}}

\renewcommand\d{{\rm d}\,}
\def\rn#1{{\bf R}^{#1}}
\def\ov#1{\overline{#1}}
\newcommand\leb{\mathcal{L}}
\newcommand\haus{\mathcal{H}}
\newcommand\Z{{\mathbb Z}}
\newcommand\N{{\mathbb N}}
\newcommand\trac{\mathop{\rm Tr}\,}

\newcommand\res{\mathop{\hbox{\vrule height 7pt width .5pt depth 0pt
\vrule height .5pt width 6pt depth 0pt}}\nolimits}
\newcommand\supp{{\rm supp}\,}  
\newcommand\Inj{{\textrm{Inj}\,}}

\newcommand\F{{\mathcal{F}}}
\newcommand\La{{\Lambda}}
\newcommand\co{{\mathcal{CO}\,}}

\newcommand\An{{\mathcal{AN}}}
\newcommand\vo{{\mathcal{V}_\infty}}
\newcommand\vd{{\mathfrak{d}\,}}

\newcommand\Is{{\mathfrak{Is}}}
\newcommand\gen{{\bf g}}

\newcommand\Ind{{\textrm{Ind}\,}}

\begin{document}

\title[The min--max construction of minimal 
surfaces]
{The min--max construction of minimal 
surfaces}

\author{Tobias H. Colding}%
\address{Courant Institute of Mathematical Sciences\\
251 Mercer Street\\ New York, NY 10012}
\author{Camillo De Lellis}%
\address{Max-Planck-Institute for Mathematics in the Sciences\\
Inselstr. 22 - 26, 04103 Leipzig / Germany}

\thanks{The first author was partially supported by NSF Grant DMS
0104453}

\email{colding@cims.nyu.edu and delellis@mis.mpg.de}

\maketitle

\centerline{Dedicated to Eugenio Calabi on occasion of his eightieth birthday} 

\section{Introduction}

In this paper we survey with complete proofs
some well--known, but hard to find, results about 
constructing closed embedded minimal surfaces in a closed 
$3$-dimensional manifold 
via min--max arguments.  This includes results 
of J. Pitts, F. Smith, and L. Simon and F. Smith.  

\begin{figure}[htbp]
\begin{center}
    \input{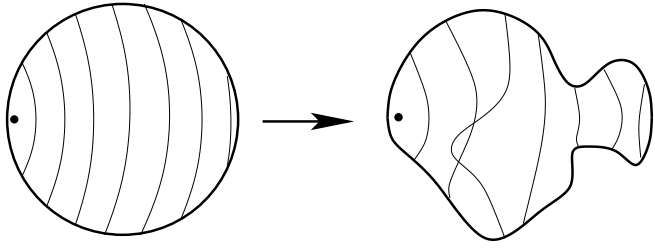}
    \caption{A 1--parameter family of curves on a $2$--sphere
    which induces a map $F:\SS^2\to \SS^2$ of degree 1.}
    \label{f:10}
\end{center}
\end{figure}

The basic idea of constructing minimal surfaces via min--max arguments
and sweep-outs goes back to Birkhoff, who used such a method to find simple 
closed
geodesics on spheres. In particular when $M^2$ is the $2$-dimensional
sphere we can find a $1$--parameter family of curves starting and ending
at a point curve in such a way that the induced map $F:\SS^2\to \SS^2$
(see Fig.~\ref{f:10}) has nonzero degree. 
Birkhoff's argument  (or the min-max
argument) allows us to conclude that $M$ has a nontrivial closed geodesic
of length less than or equal to the length of the longest curve in the
$1$-parameter family. A curve shortening argument gave that the
geodesic obtained in this way is simple.

The difficulty in generalizing this method to get
embedded minimal surfaces in $3$--manifolds is three fold. The first problem
is getting regularity
of the min--max surface obtained.
In Birkhoff's case (curves in surfaces) this was almost immediate.
The second key difficulty is
to show that the min--max surface 
is embedded. Using the technical tools of Geometric Measure
Theory (mostly the theory of varifolds), these two problems are
tackled at the same time. The third key difficulty is to get a good genus
bound for the embedded minimal surface obtained.

\subsection{The min--max construction in $3$--manifolds.}
In the following $M$ denotes a closed $3$--dimensional Riemannian manifold,
$\diff$ is the identity component of the diffeomorphism group of $M$, and 
$\Is$ is the set
of smooth isotopies.  Thus $\Is$ is the set of maps 
\ $\psi\in C^\infty ([0,1]\times M, M)$ such that
$\psi(0, \cdot)$ is the identity and $\psi(t, \cdot)\in \diff$ 
for every $t$.

We will first give a version of $1$--parameter families of surfaces
in $3$-manifolds.  
The most direct way of doing this is to let $F:[0,1]\times \Sigma\to M$ 
be a smooth map such that 
$F(t, \cdot)$ 
is an embedding of the surface $\Sigma$ 
for every $t\in [0,1]$. If we let $\Sigma_t=F (\{t\}\times \Sigma)$,
then $\{\Sigma_t\}_{t\in [0,1]}$
is a smooth 1--parameter
family of surfaces in $M$. 
This notion can be generalized in two directions.
The first one is to relax the regularity required in the 
$t$--variable:

\begin{definition}
A family 
$\{\Sigma_t\}_{t\in [0,1]}$ of surfaces of $M$ is said to be 
{\em continuous}\/ if 
\begin{itemize}
\item[(c1)] $\haus^2 (\Sigma_t)$ is a continuous function of $t$;
\item[(c2)] $\Sigma_t\to \Sigma_{t_0}$ in the
  Hausdorff topology whenever $t\to t_0$. 
\end{itemize}
\end{definition}

A second generalization allows the family of surfaces to degenerate
in finitely many points:

\begin{definition}\label{d:gensur}
A family $\{\Sigma_t\}_{t\in [0,1]}$ of subsets of $M$
is said to be a {\em generalized family}\/ of surfaces if there are a 
finite subset $T$ of $[0,1]$ and a finite set of points $P$ in $M$ such that\ 
\begin{itemize}
\item[1.] (c1) and (c2) hold;
\item[2.] $\Sigma_t$ is a surface for every $t\not \in T$;
\item[3.] For $t\in T$, $\Sigma_t$ is a surface in $M\setminus P$.
\end{itemize}
\end{definition}

Figure \ref{f:10} gives (in one dimensions less)
an example of a generalized 1--parameter family
with $T=\{0,1\}$.
To avoid confusion, families of surfaces will be denoted by 
$\{\Sigma_t\}$. Thus, when referring to a surface a subscript 
will denote a real parameter, whereas a superscript will 
denote an integer as in a sequence.

Given a generalized family $\{\Sigma_t\}$ we can generate new generalized 
families via the following procedure.
Take an arbitrary map $\psi\in C^\infty 
([0,1]\times M, M)$ such that $\psi(t, \cdot)\in \diff$ for each $t$ and
define $\{\Sigma'_t\}$ by $\Sigma'_t=\psi (t, \Sigma_t)$.\label{i:proc}
We will say that a set $\Lambda$ of generalized families is 
{\em saturated}\/ if it is closed under this operation. 

\begin{remark}\label{r:P}
For technical reasons we will require that any of the saturated sets $\Lambda$ 
that we consider has the additional property that there exists some $N=N(\Lambda)<\infty$ such that 
for any $\{\Sigma_t\}\subset \Lambda$, the set $P$ in Definition \ref{d:gensur} consists of 
at most $N$ points.  
This additional property will play a crucial role in the proof of Theorem \ref{t:SS}.
\end{remark}

Given a family $\{\Sigma_t\}\subset \Lambda$ we denote by $\F (\{\Sigma_t\})$
the area of its maximal slice and by $m_0 (\Lambda)$
the infimum of $\F$ taken over all families of $\Lambda$; that is, 
\begin{eqnarray}\label{e:min--max}
&&\F (\{\Sigma_t\}) = \max_{t\in [0,1]} \haus^2 (\Sigma_t) 
\qquad \text{ and }\\
&&m_0 (\Lambda) = \inf_{\Lambda} \F =
\inf_{\{\Sigma_t\}\in \Lambda}\, \left[ \max_{t\in [0,1]} \haus^2 
(\Sigma_t)\right]\,. 
\end{eqnarray}

\begin{figure}[htbp]
\begin{center}
    \input{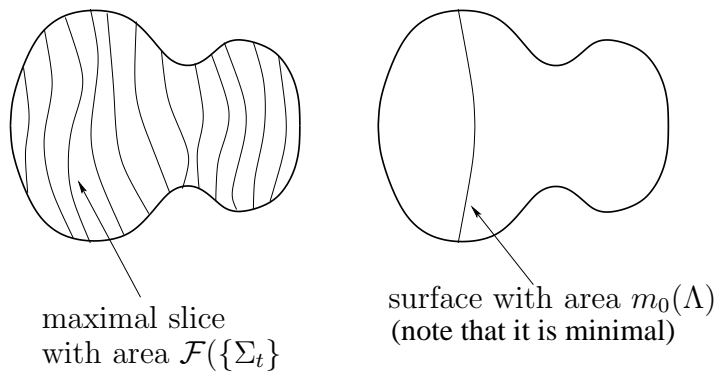}
    \caption{$\F (\{\Sigma_t\})$ and $m_0 (\Lambda)$.}
    \label{f:12}
\end{center}
\end{figure}

If $\lim_n \F (\{\Sigma_t\}^n)=m_0 (\Lambda)$, then we say that 
the sequence of generalized families of surfaces 
$\{\{\Sigma_t\}^n\}\subset \Lambda$ is a 
{\em minimizing sequence}. Assume $\{\{\Sigma_t\}^n\}$ is a minimizing 
sequence and let $\{t_n\}$ 
be a sequence of parameters. If the areas of the slices $\{\Sigma^n_{t_n}\}$ converge 
to $m_0$, i.e. if $\haus^2 (\Sigma^n_{t_n})\to m_0 (\Lambda)$, 
then we say that $\{\Sigma^n_{t_n}\}$
is a {\em min--max sequence}.

An important point in the min--max construction is to find a saturated
set $\Lambda$ of generalized families of surfaces with 
$m_0 (\Lambda)>0$. 
This can for instance be done by using the following 
elementary proposition proven in Appendix 
\ref{a:iso}; see Fig.~\ref{f:14}:

\begin{propos}\label{p:morse}
Let $M$ be a closed $3$-manifold with a Riemannian metric and let 
$\{\Sigma_t\}$ be the level sets of a Morse function. The
smallest saturated set 
$\Lambda$ containing the family $\{\Sigma_t\}$
has $m_0 (\Lambda)>0$.
\end{propos}

\begin{figure}[htbp]
\begin{center}
    \input{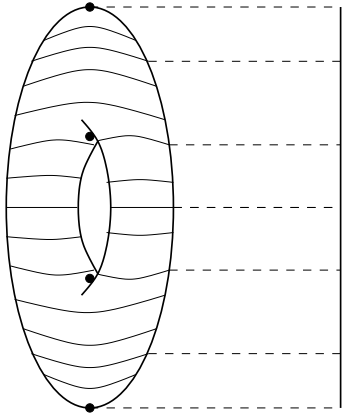}
    \caption{A sweep--out of the torus by level sets of a
      Morse function. In this case there are four degenerate slices in
      the 1--parameter family.}
    \label{f:14}
\end{center}
\end{figure}

The following example of sweep--outs of the $3$--sphere is a 
direct generalization of the families of curves on the $2$--sphere  
considered by Birkhoff:

\begin{Exa}\label{e:S2S3} Let $x_4$ be the coordinate function on the
  $3$--sphere coming from its standard embedding into $\rn{4}$. 
  By Proposition~\ref{p:morse}, for any fixed metric on $\SS^3$
  the level sets of $x_4$ generate a saturated set of generalized
  families of surfaces with $m_0>0$.
\end{Exa}

In this survey we will prove the following theorem:

\begin{theorem}\label{t:SS}[Simon--Smith] 
Let $M$ be a closed $3$-manifold with a Riemannian metric.  For any 
saturated set of generalized families of surfaces $\Lambda$, there is 
a min--max sequence obtained from $\Lambda$ and
converging in the sense of varifolds to a
smooth embedded minimal surface with area
$m_0 (\Lambda)$ (multiplicity is allowed).
\end{theorem}

An easy corollary of Proposition \ref{p:morse} and
Theorem \ref{t:SS} is the existence of a smooth embedded minimal
surface in any closed Riemannian $3$--manifold
(Pitts proved that in any closed
Riemannian manifold of dimension at most 7
there is a closed embedded 
minimal hypersurface; see theorem A and the final remark
of the introduction of \cite{P}).

For $\Lambda$ as in 
Example~\ref{e:S2S3} (where $M$ is topologically a $3$-sphere
but could have an arbitrary metric and 
the sweep--outs are by $2$-spheres) 
Simon and Smith proved that the min--max 
sequence given by Theorem~\ref{t:SS} converges to
a disjoint union of embedded minimal $2$--spheres, 
possibly with multiplicity. 
The following generalization of this result was announced by Pitts and 
Rubinstein in \cite{PR1} (in this 
theorem $\gen (\Sigma)$ is the genus of the surface $\Sigma$): 

\begin{theorem}\label{t:PR} If $\{\Sigma^k_{t_k}\}$ 
is the min--max sequence of Theorem~\ref{t:SS} and $\Sigma^\infty$
its limit, then  
\begin{equation}\label{e:PR}
\gen (\Sigma^\infty) \leq \liminf_{k\to \infty} \gen (\Sigma^k_{t_k})\,.
\end{equation}
\end{theorem}

We plan to address the proof of Theorem \ref{t:PR} elsewhere.

\part{Overview of the proof}

\section{Preliminaries}\label{s:pr} 

\subsection{Notation} We begin by fixing some notation 
which will be used throughout.  
When speaking of an isotopy $\psi$, that is, of a map $\psi: [0,1]\times
M\to M$ such that $\psi(t, \cdot)\in \diff$ for every $t$, then 
if not otherwise specified we assume that $\psi\in \Is$.  Recall 
that $\Is$ is the set of smooth isotopies that start at the identity.  
Moreover, we
say that $\psi$ {\em is supported in $U$}\/ if $\psi(t,x)=x$ for every
$(t,x)\in [0,1]\times (M\setminus U)$. 

Most places $\Gamma$ and $\Sigma$ will either denote smooth closed surfaces in $M$ 
(multiplicity allowed) or smooth surfaces in some subset $U\subset M$ with $\overline{\Sigma}\setminus
\Sigma\subset \partial U$. However, there are a few places where $\Sigma$ and $\Gamma$ denote 
surfaces which are smooth away from finitely many (singular) points.

Below is a list of our notation:

\medskip

\begin{tabular}{lll}
$T_xM$ && the tangent space of $M$ at $x$\\
$TM$ && the tangent bundle of $M$.\\
$\Inj(M)$ && the injectivity radius of $M$.\\
$\haus^2$ && the $2$--d Hausdorff 
measure in the metric\\
&& space $(M,d)$.\\
$B_\rho (x)$ && open ball\\
$\ov{B}_\rho (x)$ && closed ball\\
$\partial B_\rho (x)$ && distance sphere of radius $\rho$ in $M$.\\
$\diam (G)$ && diameter of a subset $G\subset M$.\\
$d(G_1, G_2)$ && the Hausdorff distance between the subsets\\
&& $G_1$ and $G_2$ of $M$.\\
$\D$, $\D_\rho$ && the unit disk and the disk of radius $\rho$ in $\rn{2}$.\\
$\B$, $\B_\rho$ && the unit ball and the ball of radius $\rho$ in $\rn{3}$.\\
$\exp_x$ && the exponential map in $M$ at $x\in M$.\\
$\Is (U)$ && smooth isotopies supported in $U$.\\
$G^2(U)$, $G(U)$ && grassmannian of (unoriented) 
$2$--planes\\
&& on $U\subset M$.
\end{tabular}

\medskip

\begin{tabular}{lll}
$\an (x,\tau, t)$ && the open annulus 
$B_t (x)\setminus \overline{B}_\tau (x)$.\\
$\An_r (x)$ && the set 
$\{\an (x, \tau, t) \mbox{ where $0<\tau<t<r$}\}$.\\
$C^\infty (X,Y)$ && smooth maps from $X$ to $Y$.\\
$C^\infty_c (X,Y)$ && smooth maps with compact support
from $X$\\
&& to the vector space $Y$.
\end{tabular}

\subsection{Varifolds} We will need to recall some basic facts from the 
theory of varifolds; see for instance  
chapter 4 and chapter 8 of \cite{Si} for further information.
Varifolds are a convenient way of generalizing surfaces to a category that 
has good compactness properties. An advantage of varifolds, over other
generalizations (like currents), is that they do not allow for
cancellation of mass. This last property is fundamental for the
min--max construction.

If $U$ is an open 
subset of $M$, any finite nonnegative measure on the Grassmannian of 
unoriented $2$--planes on $U$ is said to be a {\em $2$--varifold in $U$}. 
The Grassmannian 
of $2$--planes will be denoted by $G^2(U)$ and the vector space of 
$2$--varifolds is denoted by $\mathcal{V}^2 (U)$.  
With the exception of Appendix~\ref{App1}, throughout 
we will consider only $2$--varifolds; thus we drop the 2.

We endow $\mathcal{V} (U)$ with the topology of
the weak convergence in the sense of
measures, thus we say that a sequence $V^k$ of varifolds converge to
a varifold $V$ if for every function $\varphi\in C_c (G(U))$ 
$$
\lim_{k\to \infty} \int \varphi (x, \pi)\, dV^k (x, \pi)
\;=\; \int \varphi (x, \pi)\, dV (x, \pi)\, .
$$
Here $\pi$ denotes a $2$--plane of $T_x M$.
If $U'\subset U$ and $V\in \mathcal{V} (U)$, then we denote by 
$V\res U'$ the restriction of the measure $V$ to $G (U')$. Moreover, 
$\|V\|$ will be the unique measure on $U$ satisfying
$$
\int_U \varphi (x) \,d\|V\| (x)\;=\;
\int_{G(U)} \varphi (x) \,dV (x, \pi)\qquad \forall \varphi\in C_c
(U)\, .
$$
The support of $\| V\|$, denoted by $\supp (\|V\|)$, is the 
smallest closed set outside
which $\|V\|$ vanishes identically.
The number $\|V\|(U)$ will be
called the {\em mass of $V$ in $U$}. When $U$ is clear from the context, 
we say briefly the {\em mass of $V$}.

Recall also that a $2$--dimensional rectifiable set is a countable union
of closed subsets of $C^1$ surfaces (modulo sets of $\haus^2$--measure 0).
Thus, if $R\subset U$ is a $2$--dimensional rectifiable set 
and $h:R\to \rn{+}$ is a Borel function, then we can define a 
varifold $V$ by 
\begin{equation}\label{e:defvar}
\int_{G (U)} \varphi (x, \pi) \,dV (x, \pi)=
\int_R h(x) \varphi (x, T_x R) \,d\haus^2 (x)\, \quad \forall
\varphi\in C_c (G (U))\, .
\end{equation}
Here $T_x R$ denotes the tangent plane to $R$ in $x$.
If $h$ is integer--valued, then we say that $V$ is an 
{\em integer rectifiable varifold}.
If $\Sigma=\bigcup n_i \Sigma_i$, then 
by slight abuse of notation we use $\Sigma$ for the 
varifold induced by $\Sigma$ via \eqref{e:defvar}.

\subsection{Pushforward, first variation, monotonicity formula}
If $V$ is a varifold induced by a surface 
$\Sigma\subset U$ and $\psi:U\to U'$ a diffeomorphism, 
then we let $\psi_\sharp V\in \mathcal{V} (U')$ 
be the varifold induced by the 
surface $\psi (\Sigma)$. The definition of $\psi_\sharp V$ can be naturally 
extended to {\em any} $V\in \mathcal{V} (U)$ by
$$
\int \varphi(y, \sigma)\, d(\psi_\sharp V) (y, \sigma)
\;=\; \int J \psi (x, \pi)\, \varphi 
(\psi (x), d\psi_x (\pi))\, dV (x, \pi)\, ;
$$
where $J \psi (x, \pi)$ denotes the Jacobian determinant (i.e. the area
element) of the differential $d\psi_x$ restricted to the plane $\pi$;
cf. equation (39.1) of \cite{Si}.

Given a smooth vector field $\chi$, let $\psi$ be the isotopy
generated by $\chi$, i.e. with ${\textstyle \frac{\partial
    \psi}{\partial t} = \chi( \psi)}$. The   
first variation of $V$ with respect to $\chi$ is
defined as
$$
[\delta V] (\chi) \;=\; \left. \frac{d}{dt} (\|\psi (t, \cdot)_\sharp V\|)
\right|_{t=0}\, ;
$$
cf. sections 16 and 39 of \cite{Si}. When $\Sigma$ is a smooth surface 
we recover the classical definition of first variation of a
surface:
$$
[\delta \Sigma] (\chi) \;=\; \int_\Sigma {\rm div}_{\Sigma} \chi\, d\haus^2 
\;=\; \left. \frac{d}{dt} (\haus^2 (\psi(t,\Sigma)))\right|_{t=0}\, .
$$
If $[\delta V] (\chi)=0$ for every $\chi\in C^\infty_c (U,TU)$, then $V$ 
is said to be {\em stationary in $U$}. Thus stationary varifolds are 
natural generalizations of minimal surfaces.

Stationary varifolds in Euclidean spaces satisfy the monotonicity
formula (see sections 17 and 40 of \cite{Si}):
\begin{equation}\label{e:MonFor1}
\mbox{For every $x$ the function } 
f(\rho)= \frac{\|V\| (B_\rho (x))}{\pi \rho^2}
\mbox{ is non--decreasing.} 
\end{equation}
When $V$ is a stationary varifold in a Riemannian manifold a similar
formula with an error term holds. Namely, there exists a constant
$C (r)\geq 1$ such that
 \begin{equation}\label{e:MonFor}
f(s)\;\leq\; C(r) f(\rho) \qquad \mbox{whenever $0<s<\rho<r$.}
\end{equation}
Moreover, the constant $C(r)$ approaches $1$ as 
$r\downarrow 0$. This property allows us to define the
{\em density} of a stationary varifold $V$ at $x$, by
$$
\theta (x, V)\;=\; \lim_{r\downarrow 0} \frac{\|V\| (B_r (x))}{\pi r^2}.
$$  
Thus $\theta (x, V)$ corresponds to the upper
density $\theta^{*2}$ of the measure $\|V\|$ as defined in section 3
of \cite{Si}. The following theorem 
gives a useful condition for rectifiability 
in terms of density:

\begin{theorem}\label{Rectif} (Theorem 42.4 of \cite{Si}).  If $V$ is
  a stationary varifold with  
  $\theta (V,x)>0$ for $\|V\|$--a.e. $x$, then $V$ is
  rectifiable. 
\end{theorem}

\subsection{Tangent cones, Constancy Theorem} 
Tangent varifolds are the natural generalization of tangent planes for 
smooth surfaces. In order to define tangent varifolds in a $3$--dimensional
manifold we need to recall what a dilation in a manifold is. 
If $x\in M$ and $\rho<\Inj
(M)$, then the dilation around $x$ with factor $\rho$ is the map
$T^x_\rho: B_\rho (x)\to \B_1$ given by
$T^x_\rho (z)=(\exp_x^{-1} (z))/\rho$; thus if $M=\rn{3}$, then
$T^x_\rho$ is the usual dilation $y\to (y-x)/\rho$.
 
\begin{definition}\label{cones}
If $V\in \mathcal{V} (M)$, then we denote by 
$V^x_\rho$ the dilated varifold in $\mathcal{V} (\B_1)$ given by 
$V^x_\rho= (T^x_\rho)_\sharp V$.
Any limit $V'\in\mathcal{V} (\B_1)$ of a sequence $V^x_{s_n}$ 
of dilated varifolds, with $s_n \downarrow 0$, is
said to be a {\em tangent varifold at $x$}. 
The set of all tangent varifolds to $V$ at $x$ is denoted by $T(x,V)$.
\end{definition} 
It is well known that if the varifold 
$V$ is stationary, then any tangent varifold 
to $V$ is a stationary {\em Euclidean}\/ cone (see section 42 of 
\cite{Si}); that is a stationary 
varifold in $\RR^3$ which is invariant under the dilations
$y\to y/\rho$. If $V$ is also integer rectifiable and the support of $V$ 
is contained in the union of a finite number of disjoint connected surfaces 
$\Sigma_i$, i.e. $\supp (\|V\|)\subset \bigcup\, \Sigma_i$, then the
Constancy Theorem (see theorem 41.1 of \cite{Si})  gives that 
$V=\bigcup\, m_i \Sigma^i$ for some natural numbers $m_i$.

\subsection{Curvature estimates for stable minimal surfaces}
In many of the proofs we will use
Schoen's curvature estimate (see \cite{Sc} or \cite{CM2}) 
for stable minimal surfaces.  Recall that this estimate asserts that if  
$U\subset\subset M$, then there exists a universal
  constant, $C(U)$, such that for every stable minimal surface 
$\Sigma\subset U$ 
with $\partial \Sigma\subset \partial U$ and second fundamental form $A$ 
\begin{equation}\label{e:curva1}
|A|^2 (x) \; \leq\; \frac{C(U)}{d^2 (x, \partial U)}\, \qquad \forall x\in
\Sigma\, .
\end{equation}
In fact, what we will use is not the actual curvature estimate, rather 
it is the following consequence of it:
\begin{eqnarray}
&\mbox{If $\{\Sigma^n\}$ is a sequence of stable
minimal surfaces in $U$, then a}&\nonumber\\
&\mbox{subsequence converges 
to a stable minimal surface $\Sigma^\infty$}\, .&\label{e:curvaclaim}
\end{eqnarray}

\section{Overview of the proof of Theorem \ref{t:SS}}\label{s:ov}

In the following we fix a saturated set $\Lambda$ 
of generalized $1$-parameter families of surfaces and 
denote by $m_0=m_0 (\Lambda)$ the infimum
of the areas of the maximal slices in
$\Lambda$; cf. \eqref{e:min--max}.
The proof of Theorem~\ref{t:SS}, which we will outline in this section, 
follows by combining two results, Proposition~\ref{existAM} and 
Theorem~\ref{goal}. The proofs of these two results will involve all 
the material presented in Sections 3, 4, 5, and 6.  

\subsection{Stationarity} 
\begin{figure}[htbp]
\begin{center}
    \input{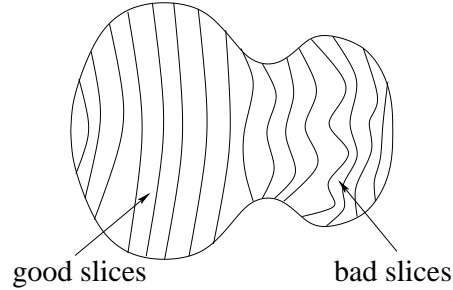}
    \caption{Slices with area close to $m_0$. The good ones are very
      near to a minimal surface of area $m_0$, whereas the bad ones
      are far from any stationary varifold.}
    \label{f:13}
\end{center}
\end{figure}

If $\{\{\Sigma_t\}^k\}\subset \Lambda$ is a 
minimizing sequence, then it is easy to show the existence 
of a min--max sequence which converge (after possibly passing to 
subsequences) to a stationary varifold. 
However, as Fig.~\ref{f:13} illustrate,
a general minimizing sequence $\{\{\Sigma_t\}^k\}$ can have
slices $\Sigma^k_{t_k}$ with area converging to $m_0$
but not ``clustering'' towards stationary 
varifolds. 

In the language introduced above, 
this means that a given minimizing sequence 
$\{\{\Sigma_t\}^k\}$ can have
min--max sequences which are not clustering to stationary varifolds.
This is a source 
of some technical problems and forces us in Section~\ref{s:stat} to show 
how to choose a ``good'' minimizing 
sequence $\{\{\Sigma_t\}^k\}$.  This is the content of the 
following proposition:

\begin{propos}\label{p:goodbis} There exists a minimizing sequence 
$\{\{\Sigma_t\}^n\}\subset \Lambda$ 
such that\ every min--max sequence $\{\Sigma^n_{t_n}\}$ 
clusters to stationary 
varifolds.
\end{propos}

A result similar to
Proposition~\ref{p:goodbis} appeared in \cite{P} (see theorem 4.3 of [P]).
The proof follows from ideas of \cite{Alm} (cf. 12.5 there).

\subsection{Almost minimizing} 
\begin{figure}[htbp]
\begin{center}
    \input{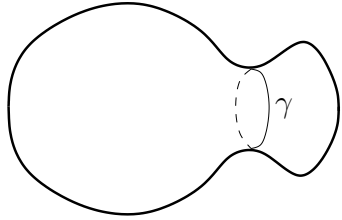}
    \caption{Curves near $\gamma$ are $\eps$--a.m.: It is 
      impossible to deform any such curve isotopically to a much smaller 
      curve without passing through a large curve.}
    \label{f:11}
\end{center}
\end{figure}

A stationary varifold can be quite far
from an embedded minimal surface. The key point for getting
regularity for varifolds produced by min--max sequences 
is the concept of ``almost minimizing surfaces'' or a.m. surfaces. 
Roughly speaking
a surface $\Sigma$ is almost minimizing if 
any path of surfaces $\{\Sigma_t\}_{t\in
  [0,1]}$ starting at $\Sigma$ and such 
that $\Sigma_1$ has small area (compared to $\Sigma$) must necessarily 
pass through a surface with large area.
That is, there must exist a $\tau\in ]0,1[$ such that\ $\Sigma_\tau$
has large area compared with $\Sigma$; see Fig.~\ref{f:11}.

The precise definition of a.m. surfaces is the following:

\begin{definition}\label{AM1} Given $\eps>0$, 
  an open set $U\subset M^3$, and a
  surface $\Sigma$, we say that $\Sigma$ is {\em $\eps$--a.m.\ in $U$}\/ if
  there {\sc does not} exist any isotopy $\psi$ supported in $U$ such that\
  \begin{eqnarray}
  &&\mbox{$\haus^2 (\psi (t,N))\leq \haus^2 (N)+\eps/8$ for all
  $t$;}\label{AM(a)}\\ 
  &&\mbox{$\haus^2 (\psi (1,N))\leq \haus^2 (N)-\eps$.}\label{AM(b)}
  \end{eqnarray}
A {\em sequence}
$\{\Sigma^n\}$ is said to be
{\em a.m.\ in $U$}\/ if each 
$\Sigma^n$ is $\eps_n$--a.m.\ in $U$ for some $\eps_n\downarrow 0$.
\end{definition}

This definition first appeared in Smith's dissertation, 
\cite{Sm}, and was inspired by a similar one
of Pitts (see the definition of almost minimizing varifolds 
in 3.1 of \cite{P}). In section 4 of his book, 
Pitts used combinatorial arguments (some of which were based on ideas
of Almgren, \cite{Alm}) to prove a general 
existence theorem for almost minimizing varifolds. The situation we deal with 
here is much simpler, due to the fact that we only consider 
$1$--parameter families of surfaces and not general multi--parameter 
families. Using a version of the
combinatorial arguments of Pitts,
we will prove in Section~\ref{s:AM} the following proposition:

\begin{propos}\label{p:existbis}
There exists a function $r:M\to \rn{+}$ and a min--max sequence
$\{\Sigma^j\}$ such that:
\begin{itemize}
\item $\{\Sigma^j\}$ is a.m. in every annulus 
$\an$ centered at $x$ and with outer radius at most $r(x)$;
\item In any such annulus, $\Sigma^j$ is smooth when $j$ is sufficiently large;
\item $\Sigma^j$ converges to a stationary varifold $V$ in $M$, as $j\uparrow \infty$.
\end{itemize}
\end{propos}

The reason why we work with annuli is two fold. The first is that 
we allow the generalized families to have slices with
point--singularities. The second is that even if any family of $\Lambda$
were made of smooth surfaces, then the combinatorial proof of Proposition
\ref{p:existbis} would give a point $x\in M$ in which we are forced to
work with annuli (cf. the proof of Proposition \ref{existAM}). 

\begin{figure}[htbp]
\begin{center}
    \input{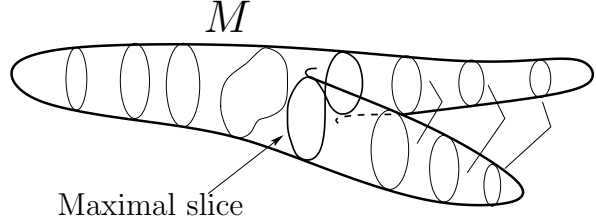}
    \caption{A sweep--out of the three--legged starfish, which can be
      realized as level--sets of a Morse function.}
    \label{f:fish}
\end{center}
\end{figure}

For a better understanding of this point 
consider the following example, due to Almgren
(\cite{Alm}, p. 15--18; see also \cite{P}, p. 20--21). 
The surface $M$ in Fig. \ref{f:fish} is
diffeomorphic to ${\bf S}^2$ and metrized as a ``three--legged
starfish''. The picture shows a  sweep--out with a unique maximal
slice, which is a geodesic figure--eight (cf. Fig. 5 of \cite{P}). 
The slices
close to the figure--eight are {\em not} almost minimizing in
balls centered at its singular point $P$. But they are almost
minimizing in every sufficiently small annulus centered at $P$.

If $\Lambda$ is the saturated set generated by the sweep--out of
Fig. \ref{f:fish}, then no min--max sequence generated by
$\Lambda$ converges to a simple closed geodesic. 
However, there are no similar examples of
sweep--outs of $3$--dimensional manifolds by $2$--dimensional objects:
The reason for this is that point--singularities of ($2$--dimensional) 
minimal surfaces are removable.

\subsection{Gluing replacements and regularity} The task  
of the last sections is to prove  that the stationary varifold $V$ 
of Proposition~\ref{p:existbis} is a smooth surface.
In Section~\ref{s:reg2} we will see 
that if $\an$ is an annulus in which $\{\Sigma^j\}$ is
a.m., then there exists
a stationary varifold $V'$, referred to as a {\em replacement}, such that\
\begin{eqnarray}
&&\mbox{$V$ and $V'$ have the same mass and $V=V'$
on $M\setminus \an$.}\label{r1}\\
&&\mbox{$V'$ is a stable minimal surface inside $\an$.}\label{r2}
\end{eqnarray}
In Lemma~\ref{rect} we use
this ``replacement property'' 
and \eqref{e:curvaclaim} to show 
that the stationary varifold $V$ of Proposition~\ref{p:existbis} 
is integer rectifiable. 
The properties of (smooth)
minimal surfaces would naturally lead to the following
unique continuation--type conjecture, cf. \cite{SW}:

\begin{Con}\label{c:unique} Let $\rho>0$ be smaller than the 
convexity radius and let $V$ and $V'$ be stationary integer 
rectifiable in $M$ with the same mass. 
If the outer radius of the annulus $\an$ is less 
than $\rho$,  $V=V'$ in $M\setminus \an$, and $V'$ is a
stable minimal surface in $\an$, then $V=V'$.
\end{Con}

An affirmative answer to Conjecture \ref{c:unique} would immediately yield 
the regularity of the stationary varifold $V$ of Proposition~\ref{p:existbis} 
in sufficiently small annuli. 
By letting the inner radius of such annuli go to zero, we would be able to
conclude that $V$ is a stable minimal surface in $B_{\rho (x)} (x)\setminus
\{x\}$, provided that $\rho (x)$ is sufficiently small. 
Hence, after showing that $x$ 
is a removable singularity we 
would get that $V$ is an embedded minimal surface.

Unluckily we are not able to argue in this way.
In fact, in 
Appendix~\ref{App1} we give an example of two distinct integer rectifiable 
$1$--varifolds $V_1$ and $V_2$ in $\rn{2}$ which have 
the same mass and coincide outside a disk. This example does 
not disprove Conjecture~\ref{c:unique}; because, besides the
dimensional difference, in the disk 
where $V_1\neq V_2$, both the varifolds are singular. 
It does however show that a proof of Conjecture~\ref{c:unique} could 
be rather delicate.

In \cite{Sm} this problem of unique continuation
was overcome by showing that for $V$ as in 
Proposition~\ref{p:existbis}, 
one can construct ``secondary'' replacements $V''$ also for the 
replacements $V'$. This idea goes back to \cite{P}. 
In Section~\ref{s:reg1} we follow \cite{Sm} and show that if we can replace
sufficiently many times, then $V$ is regular
(cf. Definition~\ref{goodrepl}
and Proposition~\ref{reg} for the precise statement). 

\subsection{Replacements} As discussed in the previous subsection, 
to prove the regularity of $V$, we need 
to construct (sufficiently many) replacements.
This task is accomplished in two steps in Section~\ref{s:reg2}. 

{\bf Step 1}: Fix an annulus $\an$ in $M$ in which $\Sigma^k$ is
$\eps_k$--a.m.\ In this annulus we deform $\Sigma^k$ into a further
sequence of surfaces $\{\Sigma^{k,l}\}^l$ with the following properties:
\begin{itemize}
\item $\Sigma^{k,l}$ is the image of $\Sigma^k$ under some isotopy
  $\psi$ which satisfies \eqref{AM(a)} (with
  $\eps=\eps_k$ and $U=\an$);
\item If we denote by $\mathcal{S}^k$ the family of all such isotopies,
  then
\begin{equation}\label{e:minp}
\lim_{l\to \infty} \haus^2 (\Sigma^{k,l})=\inf_{\psi\in \mathcal{S}^k}
\haus^2 (\psi (1, \Sigma^k))\,.
\end{equation}
\end{itemize}
After possibly passing to a subsequence, then 
$\Sigma^{k,l}\to V^k$ and $V^k\to V'$,
where $V^k$ is a varifold which is stationary {\em in $\an$}.
By the a.m. property of $V$, it follows that  
$V'$ is stationary {\em in all of $M$}\/ and
satisfies \eqref{r1}.

The second step is to prove that 
$V^k$ is a (smooth) stable minimal surface in $\an$. Thus,
\eqref{e:curvaclaim} will give that also
$V'$ is a stable minimal surface in $\an$. After checking some details we  
show that $V$ meets the technical requirements of
Proposition~\ref{reg}.

{\bf Step 2}: It remains to prove that $V^k$ is a stable minimal
surface. Stability is a trivial consequence of \eqref{e:minp}. For the
regularity we use again Proposition~\ref{reg}. The key is proving the
following property: 
\begin{itemize}
\item[(P)] If $B\subset \an$ is a sufficiently small ball and
$l$ is a sufficiently large number, then {\em any} $\psi\in \Is (B)$
with $\haus^2 (\psi (1, \Sigma^{k,l}))\leq 
\haus^2 (1, \Sigma^k)$ can be replaced by a $\Psi\in \Is (\an)$ with 
\begin{equation}\label{e:P}
\Psi(1, \cdot) = \psi (1, \cdot)\quad
\mbox{and} \quad \haus^2 (\Psi (t, \Sigma^{k,l}))\leq 
\haus^2 (\Sigma^{k,l}) +
  \eps_k/8\, \quad \mbox{for all $t$}.
\end{equation}
\end{itemize}

We will now discuss how (P) gives the regularity of $V^k$.

Fix a sufficiently small ball $B$ and a large number $l$ so 
that the property (P) above holds.
Take a sequence of surfaces $\Gamma^j=\Sigma^{k,l,j}$ which are 
isotopic to $\Sigma^{k,l}$ in $B$ and
such that\ $\haus^2 (\Gamma^j)$ converges to 
$$
\inf_{\psi\in \Is (B)} \haus^2 (1, \psi(\Sigma^{k,l}))\, .
$$ 
By a result of Meeks--Simon--Yau, \cite{MSY}, $\Gamma^j$ converges 
to a varifold $V^{k,l}$ which is a 
stable minimal surface in $B$. Thus, by \eqref{e:curvaclaim}, 
the sequence of varifolds $\{V^{k,l}\}^l$ converges to a varifold $W^k$ which
is a stable minimal surface in $B$. 
The property (P) is used to show that, for $j$ and $l$ sufficiently
large, $\Sigma^{k,l,j}$ is 
a good competitor with respect to the 
$\eps_k$--a.m. property of $\Sigma^k$. This is then used to show 
that $W^k$ is a replacement for
$V^k$ in $B$. Again it is only a technical step to check that we 
can apply Proposition~\ref{reg}, and hence get that $V^k$ 
is a stable minimal surface in $\an$. 

\part{Proof of Theorem \ref{t:SS}}

\section{Limits of suitable min--max sequences are stationary}\label{s:stat}

This section is devoted to the proof of Proposition~\ref{p:goodbis}.  
For simplicity we 
metrize the weak topology on the space of varifolds and restate 
Proposition~\ref{p:goodbis} 
using this metric.

Denote by $X$ the set of
varifolds $V\in \mathcal{V} (M)$ with mass bounded by $4\,m_0$, i.e., with 
$\|V\| (M)\leq 4m_0$. Endow $X$ with the weak$^*$ topology and let $\vo$ be 
the set of stationary varifolds contained in $X$. Clearly, $\vo$ is a
closed subset of $X$. Moreover, by standard general topology theorems, 
$X$ is compact and metrizable.
Fix one such metric and denote it by $\vd$. The
ball of radius $r$ and center $V$ in  
this metric will be denoted by $U_r (V)$.

\begin{propos}\label{good} There exists a minimizing sequence 
$\{\{\Sigma_t\}^n\}\subset \Lambda$
  such that, if $\{\Sigma^n_{t_n}\}$ is a min--max sequence, then 
  $\vd (\Sigma^n_{t_n}, \vo)\to 0$.
\end{propos}
\begin{proof} The key idea of the proof is building a continuous
  map $\Psi:X\to \Is$ such that :
\begin{itemize}
\item  If $V$ is stationary, then $\Psi_V$ is the trivial isotopy;
\item If $V$ is not stationary, then $\Psi_V$
decreases the mass of $V$.
\end{itemize}
Since each $\Psi_V$ is an isotopy, and thus is itself a map from 
$[0,1]\times M\to M$, to avoid confusion we use the subscript $V$
to denote the dependence on the varifold $V$.  
The map $\Psi$ will be used to deform a minimizing sequence
$\{\{\Sigma_t\}^n\}\subset \Lambda$ into another minimizing 
sequence $\{\{\Gamma_t\}^n\}$ such that :
\begin{itemize}
\item[] For every $\eps>0$, there exist $\delta>0$ and $N\in \N$ 
such that 
\begin{equation}\label{condS}
\mbox{if}\quad \left\{
\begin{array}{ll}
& n>N\\
\mbox{and} & \haus^2 (\Gamma^n_{t_n})>m_0-\delta
\end{array}\right\}\, ,\quad
\mbox{then} \quad \vd (\Gamma^n_{t_n}, \vo)< \eps.
\end{equation}
\end{itemize}
Such a $\{\{\Gamma_t\}^n\}$ would satisfy the requirement of the proposition.

The map $\Psi_V$ should be thought of as a natural ``shortening process'' 
of varifolds which are not stationary. If the mass (considered as a 
functional on the 
space of varifolds) were smoother, then a gradient flow would provide a 
natural shortening process like $\Psi_V$. 
However, this is not the case; even if we start with 
smooth initial datum, in very short time
the motion by mean curvature, i.e.\ the gradient flow of the area functional
on smooth submanifolds, gives surfaces which are not isotopic to 
the initial one.

\medskip

\noindent
{\bf Step 1: A map from $X$ to the space of vector fields.}

The isotopies $\Psi_V$ will be generated as 1--parameter families of
diffeomorphisms satisfying certain ODE's. In this step we associate to
any $V$ a suitable vector field, which in Step 2 will be used to 
construct $\Psi_V$.

For $k\in\Z$ define the annular neighborhood of $\vo$ 
$$
\mathcal{V}_k\;=\;\left\{V\in X| 2^{-k+1}\geq\vd (V, \vo)\geq 
2^{-k-1}\right\} \,.
$$
There exists a positive constant $c(k)$ depending on $k$ such that 
to every $V\in \mathcal{V}_k$ we can associate a smooth vector field
$\chi_V$ with 
$$
\|\chi_V\|_\infty\leq 1 \qquad \mbox{and} \qquad \delta V (\chi_V)\leq -
c(k).
$$
Our next task is choosing $\chi_V$ with continuous dependence on $V$.
Note that for every $V$ there is some radius $r$ such that
$\delta W (\chi_V)\leq - c(k)/2$ for every $W\in U_r (V)$.
Hence, for any $k$ we can find balls $\{U^k_i\}_{i=1, \ldots, N(k)}$ and 
vector fields $\chi_i^k$ such that :
\begin{eqnarray}
&&\mbox{The balls $\tilde{U}^k_i$ 
concentric to $U^k_i$ with half the radii
  cover $\mathcal{V}_k$;}\quad\label{primaa}\\
&&\mbox{If $W \in U^k_i$, then $\delta W (\chi)\leq -
  c(k)/2$;}\label{primab}\\
&&\mbox{The balls $U^k_i$ are disjoint from 
$\mathcal{V}_j$ if $|j-k|\geq 2$.}\label{primac}
\end{eqnarray}
Hence, $\{U^k_i\}_{k,i}$ is a locally finite covering of
$X\setminus \vo$. To this family we can subordinate a
continuous partition of unit $\varphi^k_i$. Thus we 
set $H_V= \sum_{i,k} \varphi^k_i (V)\chi^k_i$. The
map $H:X\to C^\infty (M, TM)$ which to every $V$ associates
$H_V$ is continuous. Moreover, $\| H_V\|_\infty\leq 1$ for
every $V$. 

\medskip

\noindent
{\bf Step 2: A map from $X$ to the space of isotopies.}

For $V\in \mathcal{V}_k$ we let $r(V)$ be the radius of the smaller
ball $\tilde{U}^j_i$ which contains it. We find that $r(V)> r(k)>0$, where
$r(k)$ only depends on $k$. Moreover, by \eqref{primab} and
\eqref{primac}, for every
$W$ contained in the ball $U_{r(V)} (V)$ we have that
$$
\delta W (H_V)\;\leq\; - \frac{1}{2}\min \{ c(k-1), c(k), c(k+1)\}.
$$
Summarizing there are two continuous functions $g:\rn{+}\to \rn{+}$ and 
$r:\rn{+}\to \rn{+}$ such that 
\begin{equation}\label{lowering}
\delta W (H_V)\;\leq\; - g(\vd (V, \vo)) \qquad 
\mbox{if} \qquad \vd (W, V)\leq r(\vd (V, \vo)).
\end{equation}
Now for every $V$ construct the 1--parameter family of
diffeomorphisms
$$
\Phi_V:[0, +\infty)\times M\to M\qquad \mbox{with} \qquad
\frac{\partial \Phi_V (t,x)}{\partial t}\;=\; H_V (\Phi_V (t,x)).
$$
For each $t$ and $V$, we denote by $\Phi_V (t, \cdot)$ the
corresponding diffeomorphism of $M$. 
We claim that there are continuous
functions $T: \rn{+} \to [0,1]$ and $G: \rn{+} \to \rn{+}$ such that  
\begin{itemize}
\item[-] If $\gamma=\vd (V, \vo)>0$
  and we transform $V$ into $V'$ via the diffeomorphism $\Phi_V
  (T(\gamma), \cdot)$, then $\|V'\| (M)\leq \|V\| (M)- G(\delta)$;
\item[-] $G(s)$ and $T(s)$ both converge to $0$ as $s\downarrow 0$.
\end{itemize}
Indeed fix $V$. For every $r>0$ there is a
$T>0$ such that the curve of varifolds 
$$
\{V(t)=(\Phi_V (t, \cdot))_\sharp V\, ,\quad t\in [0,T]\}
$$
stays in $U_r (V)$. Thus
\begin{eqnarray*}
\|V (T)\| (M)- \|V\|(M)&=& \|V(T)\| (M)- \|V(0)\| (M)\\ 
&\leq& \int_0^T [\delta V (t)]
(H_V)\, dt\, , 
\end{eqnarray*}
and therefore if we choose $r= r(\vd (V, \vo))$ as in (\ref{lowering}), then
we get the bound
$$
\|\Gamma (T)\| (M) - \|V\| (M)\; \leq \; - T g(\vd (V, \vo)).
$$
Using a procedure similar to that of Step 1 we can
choose $T$ depending continuously on $V$. It is then trivial to see that
we can in fact choose $T$ so that at the same time it is continuous and 
depends only on $\vd (V, \vo)$.

\medskip

\noindent
{\bf Step 3: Constructing the competitor and the conclusion.}

For each $V$, set $\gamma=\vd (V, \vo)$ and 
$$
\Psi_V (t, \cdot )\;=\; \Phi_V ([T(\gamma)]\, t, \cdot)
\qquad \mbox{for $t\in [0,1]$.}
$$
$\Psi_V$ is a ``normalization'' of $\Phi_V$. 
From Step 2 we know that there is a continuous function
$L:\rn{}\to \rn{}$ such that
\begin{itemize}
\item[-] $L$ is strictly increasing and $L(0)=0$;
\item[-] $\Psi_V (1, \cdot)$ deforms $V$ into a varifold $V'$
with $\|V'\|\leq \|V\| - L(\gamma)$.
\end{itemize}
Choose a sequence of families $\{\{\Sigma_t\}^n\}\subset \Lambda$ with $\F
(\{\Sigma_t\}^n)\leq m_0 +1/n$ and define $\{\Gamma_t\}^n$ by
\begin{equation}
\Gamma^n_t = \Psi_{\Sigma^n_t} (1, \Sigma^n_t) 
\qquad \mbox{for all $t\in [0,1]$ and all
  $n\in\N$}
\end{equation}
Thus
\begin{equation}\label{replacement}
\haus^2 (\Gamma^n_t)\leq \haus^2 (\Sigma^n_t) - L(\vd (\Sigma^n_t, \vo)).
\end{equation}
Note that $\{\Gamma_t\}^n$ does not necessarily belong to $\Lambda$, 
since the families of diffeomorphisms
$\psi_t (\cdot) =\Psi_{\Sigma^n_t} (1, \cdot)$ may not depend smoothly 
on $t$. In order to 
overcome this technical obstruction fix $n$ and 
note that $\Psi_t=\Psi_{\Sigma^n_t}$ is the 1--parameter family of 
isotopies generated by the $1$--parameter family of vector fields 
$h_t= T (\Sigma^n_t) H_{\Sigma^n_t}$. Think of $h$ as 
a continuous map 
$$
h:[0,1] \to C^\infty (M, TM) \quad \mbox{with the topology of 
$C^k$ seminorms}. 
$$
Thus $h$ can be approximated by a {\em smooth}\/ map 
$\tilde{h}:[0,1]\to C^\infty (M, TM)$.
Consider the {\em smooth}\/
$1$--parameter family of isotopies $\tilde{\Psi}_t$ generated by the 
vector fields
$\tilde{h}_t$ and the family of surfaces $\{\Gamma_t\}^n$ given 
by $\Gamma^n_t= 
\tilde{\Psi}_t (1, \Sigma^n_t)$. If $\sup_t \|h_t-\tilde{h}_t\|_{C^1}$ 
is sufficiently small, then we 
easily get (by the same calculations of the previous steps)
\begin{equation}\label{e:35}
\haus^2 (\Gamma^n_t)\leq \haus^2 (\Sigma^n_t) - 
L(\vd (\Sigma^n_t, \vo))/2.
\end{equation}
Moreover, since $\tilde{\Psi}_t (1, \cdot)$ is a smooth map, 
this new family belongs to $\Lambda$.

Clearly $\{\{\Gamma_t\}^n\}$ is a minimizing sequence. We next show 
that $\{\{\Gamma_t\}^n\}$
satisfies \eqref{condS}.
Note first that the construction yields a 
continuous and increasing function $\lambda:\rn{+}\to \rn{+}$ such that 
\begin{equation}\label{continuity}
\mbox{$\lambda (0)=0$}\qquad \mbox{and}\qquad \vd (\Sigma^n_t, \vo)
\geq \lambda (\vd (\Gamma^n_t, \vo))\,.
\end{equation}
Fix $\eps>0$ and choose $\delta>0$, $N\in \N$ such
that $L(\lambda (\eps))/2 -\delta > 1/N$.  
We claim that \eqref{condS} is satisfied with this choice.  Suppose not; 
then there are
$n>N$ and $t$ such that $\haus^2 (\Gamma^n_t)> m_0- \delta$
and
$\vd (\Gamma^n_t, \vo)> \eps$. Hence, 
from \eqref{e:35} and \eqref{continuity} we get
$$
\haus^2 (\Sigma^n_t)\;\geq\;\haus^2 (\Gamma^n_t) + \frac{L(\lambda
(\eps))}{2} -\delta \; >\; m_0 + \frac{1}{N}\geq m_0 + \frac{1}{n}\,.
$$
This contradicts the assumption that $\F (\{\Sigma_t\}^n)\leq m_0+1/n$. 
Thus \eqref{condS} holds and the proof is completed.
\end{proof}

\section{Almost minimizing min--max sequences}\label{s:AM}

As above, $\Lambda$ is a fixed saturated set of $1$--parameter
families $\{\Sigma_t\}$ in $M$. In the previous section we showed that 
there exists a family $\{\Sigma_t\}$ such that every min--max sequence
is clustering towards stationary varifolds. We will
now prove that one of these min--max sequences is a.m. in sufficiently
many annuli. 

\begin{propos}\label{existAM}
There exists a function $r:M\to \rn{+}$ and a min--max 
sequence $\{\Sigma^j\}$ such that:
\begin{eqnarray}
&&\mbox{$\{\Sigma^j\}$ is a.m.\ in every $\an\in \An_{r(x)} (x)$, 
for all $x\in M$}\, .\label{primai}\\
&&\mbox{In every such $\an$, $\Sigma^j$ is a smooth surface}\nonumber\\ 
&&\mbox{when $j$ is sufficiently large}\, .\label{primaiii}\\
&&\mbox{$\Sigma^j$ converges to a stationary 
varifold $V$ as $j\uparrow \infty$.}\label{primaii}
\end{eqnarray}
\end{propos}

We first fix some notation.

\begin{definition}\label{AM2}
  Given a pair of open sets $(U^1,U^2)$ we say that a surface $\Sigma$
  is {\em $\eps$--a.m.\ in $(U^1,U^2)$}\/ if it is 
  $\eps$--a.m.\
  in at least one of the two open sets. 
We denote by $\co$ the set of pairs $(U^1, U^2)$ of open sets with 
$$
\d (U^1, U^2)\geq 2\min \{\diam (U^1), \diam (U^2)\}\, .
$$
\end{definition}

Proposition \ref{existAM} will be an easy corollary of the following:

\begin{propos}\label{AMinpairs}
There is a min--max sequence $\{\Sigma^L\}=\{\Sigma^{n(L)}_{t_{n(L)}}\}$ 
which converges 
to a stationary varifold and such that 
\begin{equation}\label{e:ii}
\mbox{each $\Sigma^L$ is $1/L$--a.m.\ in every $(U^1, U^2)\in \co$.}
\end{equation} 
\end{propos}

Note that the $\Sigma^L$'s in the previous proposition may be degenerate 
slices (that is, 
they may have a finite number of singular points). 
The key point for proving Proposition~\ref{AMinpairs} is the following
obvious lemma:

\begin{lemma}\label{pairs2}
If $(U^1, U^2)$ and $(V^1, V^2)\in \co$, then there are $i,j\in \{1,2\}$
with $\d (U^i, V^j)>0$.
\end{lemma}

Before giving a rigorous proof of
Proposition~\ref{AMinpairs} we will explain the ideas behind it. 

\subsection{Outline of the proof of Proposition~\ref{AMinpairs}}
First of all note that if
a slice $\Sigma^n_{t_0}$ 
is not $\eps$--a.m.\ in 
a given open set $U$, then we can decrease its area by an isotopy 
$\psi$ satisfying \eqref{AM(a)} and \eqref{AM(b)}. Now fix an open interval 
$I$ around $t_0$ and choose a smooth bump function 
$\varphi\in C^\infty_c (I, [0,1])$ 
with $\varphi (t_0)=1$. Define $\{\Gamma_t\}^n$ by
$$
\Gamma^n_t \;=\; \psi (\varphi (t), \Sigma^n_t)\,.
$$
If the interval $I$ is sufficiently small, then by \eqref{AM(a)},
for any $t\in I$,
the area of $\Gamma^n_t$ will not be much larger than 
the area of $\Sigma^n_t$.
Moreover, for $t$ very close to $t_0$ 
(say, in a smaller interval $J\subset I$) 
the area of $\Gamma^n_t$ 
will be much less than the area of $\Sigma^n_t$.

\begin{figure}[htbp]
\begin{center}
    \input{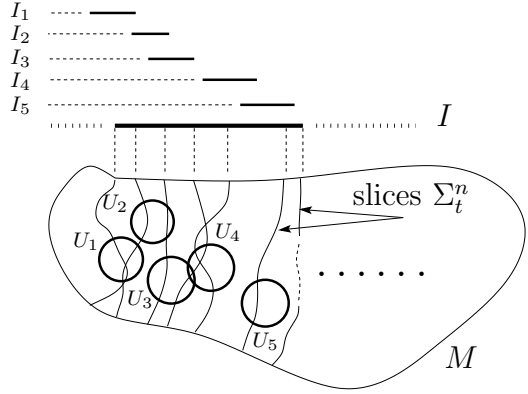}
    \caption{The covering $I_j$ and the sets $U_j$. No point of $I$ is 
contained in more than two $I_j$'s.  The intersection 
    $U_j\cap U_k=\emptyset$ if $I_j$ and $I_k$ overlap.}
    \label{f:1}
\end{center}
\end{figure}

We will show Proposition~\ref{AMinpairs} by arguing by
contradiction. So suppose that the proposition fails;
we will construct a better competitor $\{\{\Gamma_t\}^n\}$. Here 
the pairs $\co$ will play a crucial role. Indeed
when the area of $\Sigma^n_t$ is sufficiently 
large (i.e.\ close to $m_0$),  
we can find {\em two}\/ disjoint open sets $U_1$ and $U_2$ in which 
$\Sigma^n_t$ is not almost minimizing.
Consider the set $K_n\subset [0,1]$ of slices with sufficiently 
large area. Using Lemma~\ref{pairs2} (and some elementary 
considerations),
we find a finite family of intervals $I_j$, open sets $U_j$, and isotopies 
$\psi_j:I_j\times M \to M$ satisfying the 
following conditions; see Fig.~\ref{f:1}:
\begin{eqnarray}
&&\mbox{$\psi_j$ is supported in $U_j$ and is the identity at the ends 
of $I_j$.}\qquad\label{prima1}\\ 
&&\mbox{If $I_j\cap I_k\neq \emptyset$, then $U_j\cap U_k=\emptyset$.}
\label{prima2}\\
&&\mbox{No point of $[0,1]$ belong to more than two $I_j$'s.}
\label{prima3}\\
&&\mbox{$\haus^2(\psi_j (t, \Sigma^n_t))$ is never much larger than 
$\haus^2(\Sigma^n_t$).}\label{prima4}\\
&&\mbox{For every $t\in K_n$, there is $j$ s.t. 
$\haus^2(\psi_j (t, \Sigma^n_t))$ 
is much}\nonumber\\
&&\mbox{smaller than $\haus^2(\Sigma^n_t)$.}\label{prima5}
\end{eqnarray}

Conditions \eqref{prima1} and \eqref{prima2} 
allow us to ``glue'' the $\psi_j$'s in a unique 
$\psi\in \Is$ such that $\psi=\psi_j$ on $I_j\times U_j$.
The family $\{\Gamma_t\}^n$ given by $\Gamma^n_t=\psi (t, \Sigma^n_t)$
is our competitor. Indeed for every $t$,
there are at most two $\psi_j$'s which change $\Sigma^n_t$. If $t\not
\in K_n$, then none of them
increases the area of $\Sigma^n_t$ too much. Whereas, if $t\in K_n$,
then one $\psi_j$ decreases 
the area of $\Sigma^n_t$ a definite amount, and 
the other increases the area of $\Sigma^n_t$ a small amount. 
Thus, the area of the ``small--area'' slices are not increased much and 
the area of ``large--area'' slices are decreased. This yields
that $\F (\{\Gamma_t\}^n)-
\F (\{\Sigma_t\}^n)<0$. We will now give a rigorous bound for this
(negative) difference.

\subsection{Proof of Proposition~\ref{AMinpairs}}

\begin{proof}[Proof of Proposition~\ref{AMinpairs}]  
We choose 
$\{\{\Sigma_t\}^n\}\subset \La$ such that $\F
  (\{\Sigma_t\}^n)< m_0+1/n$ and satisfying the requirements of 
Proposition~\ref{good}. 
Fix $L\in \N$. To prove the proposition we claim
  there exist $n>L$ and $t_n\in [0,1]$ such that $\Sigma^n=\Sigma^n_{t_n}$ 
  satisfies \eqref{e:ii} and $\haus^2 (\Sigma^n)\geq
  m_0-1/L$. We define the sets
$$
K_n\;=\;\left\{t\in [0,1] : \haus^2 (\Sigma^n_t)\geq
  m_0-\frac{1}{L}\right\}
$$
and argue by contradiction.
Suppose not; then for every 
$t\in K_n$ there exists a pair of open subsets $(U_t^1 , U_t^2)$ such that
$\Sigma^n_t$ is not $1/L$--a.m.\ in either of them. So for
every $t\in K_n$ there exists isotopies $\psi_t^i$ such that
\begin{itemize}
\item[(1)] $\psi_t^i$ is supported on $U_t^i$;
\item[(2)] $\haus^2(\psi_t^i(1, \Sigma^n_t))\leq \haus^2
  (\Sigma^n_t)-1/L$;
\item[(3)] $\haus^2(\psi_t^i(\tau, \Sigma^n_t))\leq \haus^2 
(\Sigma^n_t)+1/(8L)$ 
for every $\tau\in [0,1]$.
\end{itemize}
In the following we fix $n$ and drop the subscript from $K_n$.
Since $\{\Sigma^n_t\}$ is continuous in $t$, 
if $t\in K$ and $|s-t|$ is sufficiently small,
then
\begin{itemize}
\item[(2')] $\haus^2(\psi_t^i(1, \Sigma^n_s))\leq \haus^2
  (\Sigma^n_s)-1/(2L)$;
\item[(3')] $\haus^2(\psi_t^i(\tau, \Sigma^n_s))\leq \haus^2 (\Sigma^n_s)
+1/(4L)$ for every $\tau\in [0,1]$.
\end{itemize}
By compactness we can cover $K$ with a finite
number of intervals satisfying (2') and (3').
This covering $\{I_k\}$ can be chosen so 
that $I_k$ overlaps only with $I_{k-1}$ and $I_{k-2}$.
Summarizing we can find 
\begin{eqnarray*}
\mbox{closed intervals} && I_1, \ldots I_r\\ 
\mbox{pairs of open sets} && (U^1_1,U^2_1), \ldots, (U^1_r, U^2_r) \in \co\\
\mbox{and pairs of isotopies} && (\psi^1_1, \psi^2_1), \ldots, (\psi^1_r,
\psi^2_r) 
\end{eqnarray*}
such that
\begin{itemize}
\item[(A)] the interiors of $I_j$ cover $K$ and $I_j\cap I_k=\emptyset$ if
  $|k-j|\geq 2$;
\item[(B)] $\psi^i_j$ is supported in $U^i_j$;
\item[(C)] $\haus^2(\psi^i_j(1, \Sigma^n_s))\leq \haus^2
  (\Sigma^n_s)-1/(2L)$\, $\forall s\in I_j$;
\item[(D)] $\haus^2(\psi^i_j(\tau, \Sigma^n_s))\leq \haus^2
  (\Sigma^n_s)+1/(4L)$\, $\forall s\in I_j$ and $\tau\in [0,1]$.
\end{itemize}
In Step 1 we refine this covering. In Step 2 we use the refined
covering to construct a competitor
$\{\Gamma_t\}^n\in \La$ with
\begin{equation}\label{e:desiredbound}
\F (\{\Gamma_t\}^n)\leq \F (\{\Sigma_t\}^n)-1/(2L)\, .
\end{equation} 
The arbitrariness of $n$ will give that
$\liminf_n \F (\{\Gamma_t\}^n)< m_0$. This is the desired
contradiction which yields the proposition.

\medskip

\noindent
{\bf Step 1: Refinement of the covering.}

First we want to find 
\begin{eqnarray*}
\mbox{a covering $\{J_1, \ldots, J_R\}$} && \mbox{which is a 
refinement of $\{I_1,
  \ldots I_r\}$}\, ,\\
\mbox{open sets $V_1, \ldots V_R$} && \mbox{among $\{U^i_j\}$}\, ,\\
\mbox{and isotopies $\varphi_1, \ldots, \varphi_R$} && \mbox{among
  $\{\psi^i_j\}$}\, ,
\end{eqnarray*}
such that: 
\begin{itemize}
\item[(A1)] The interiors of $J_i$ cover $K$ and $J_i\cap J_k=\emptyset$
  for $|k-i|\geq 2$;
\item[(A2)] If $J_i\cap J_k\neq \emptyset$, then $d(V_i, V_k)>0$;
\item[(B')] $\varphi_i$ is supported in $V_i$;
\item[(C')] $\haus^2(\varphi_i(1, \Sigma^n_s))\leq \haus^2
  (\Sigma^n_s)-1/(2L)$\, $\forall s\in J_i$;
\item[(D')] $\haus^2(\varphi_i(\tau, \Sigma^n_s))\leq \haus^2
  (\Sigma^n_s)+1/(4L)$\, $\forall s\in J_i$ and $\tau\in [0,1]$.
\end{itemize}
We start by setting $J_1=I_1$ and we distinguish two cases.
\begin{itemize}
\item[-] {\bf Case a1}: $I_1\cap I_2=\emptyset$; we set
$V_1=U^1_1$, and $\varphi_1=\psi^1_1$.
\item[-] {\bf Case a2}: $I_1\cap I_2\neq \emptyset$; by Lemma~\ref{pairs2} 
we can choose $i,k\in \{1,2\}$ such that
$\d (U^i_1, U^k_2)>0$ and we set $V_1= U^i_1$, $\varphi_1=\psi^i_1$. 
\end{itemize}

\noindent We now come to the choice of $J_3$. If we come from case
a1 then:
\begin{itemize}
\item[-] {\bf Case b1:} We make our choice as above
  replacing $I_1$ and $I_2$ with $I_2$ and $I_3$; 
\end{itemize}
If we come from case a2, then we let $i$ and $k$
be as above and we further distinguish two cases.
\begin{itemize}
\item[-] {\bf Case b21}: $I_2\cap I_3=\emptyset$; we define $J_2= I_2$,
  $V_2= U^k_2$, $\varphi_2=\psi^k_2$.
\item[-] {\bf Case b22}: $I_2\cap I_3\neq \emptyset$; by Lemma~\ref{pairs2} 
  there exist $l,m\in
  \{1,2\}$ such that $d(U^l_2, U^m_3)>0$. If $l=k$, then we define
  $J_2=I_2$, $V_2= U^k_2$, $\varphi_2=\psi^k_2$. Otherwise we choose two
  closed intervals $J_2, J_3\subset I_2$ such that 

\begin{itemize}
\item their interiors cover the interior of $I_2$,
\item $J_2$ does not overlap with any $I_h$ for $h\neq 1,2$,
\item $J_3$ does not overlap with any $I_h$ for $h\neq 2,3$.
\end{itemize}
Thus we set $V_2=U_2^k$, $\varphi_2=\psi_2^k$, and $V_3=U_2^l$,
$\varphi_3= \psi_2^l$.
\end{itemize}
An inductive argument using this procedure gives the desired
covering. Note that the cardinality of $\{J_1, \ldots, J_R\}$ is
at most $2r-1$.

\medskip

\noindent
{\bf Step 2: Construction.}

Choose $C^\infty$ functions $\eta_i$ on $\RR$ taking values in
$[0,1]$, supported in $J_i$, and such that for every $s\in K$, there
exists $\eta_i$ with $\eta_i(s)=1$. 
Fix $t\in [0,1]$ and consider the set $\Ind_t\subset \N$
of all $i$ containing $t$; thus $\Ind_t$ consists of at most two elements. 
Define 
subsets of $M$ by 
\begin{equation}\label{definition}
\Gamma^n_t\;=\;\left\{
\begin{array}{ll}
\varphi_i(\eta_i (t),\Sigma^n_t) & \mbox{in the open sets $V^i$, 
$i\in\Ind_t$\, ,}\\
\Sigma^n_t & \mbox{outside}.
\end{array}\right.
\end{equation}
In view of (A1), (A2) and (B'), then $\{\Gamma_t\}^n$ is well defined and
belongs to $\Lambda$. 

\medskip

\noindent
{\bf Step 3: The contradiction.}

We now want to bound the energy $\F (\{\Gamma_t\}^n)$ and hence we
have to estimate $\haus^2 (\Gamma^n_t)$. Note that
by (A1) every $\Ind_t$ consists of at most two integers. 
Assume for the sake of argument 
that $\Ind_t$ consists of {\em exactly} two integers.
From the construction, there exist $s_i, s_k\in [0,1]$ such that
$\Gamma^n_t$ is obtained from $\Sigma^n_t$ via the diffeomorphisms
$\varphi_i (s_i, \cdot)$, $\varphi_k (s_k, \cdot)$.
By (A2) these diffeomorphisms are supported on disjoint
sets. Thus if $t\not \in K$, then (D') gives
$$
\haus^2 (\Gamma^n_t)\;\leq\; \haus^2 (\Sigma^n_t) + \frac{2}{4L}
\;\leq\; m_0 -\frac{1}{2L}.
$$
If $t\in K$, then at least one of $s_i, s_k$ is equal to 1. Hence (C) and
(D) give
$$
\haus^2 (\Gamma^n_t)\;\leq\; \haus^2 (\Sigma^n_t) - \frac{1}{L} +\frac{1}{4L}
\;\leq\; \F (\{\Sigma^n_t\})-\frac{3}{4L}\, .
$$
Therefore $\F (\{\Gamma_t\}^n)\leq \F (\{\Sigma_t\}^n) - 1/(2L)$. This
is the desired bound \eqref{e:desiredbound}. 
\end{proof}

We now come to Proposition~\ref{existAM}. 

\begin{proof}[Proof of Proposition~\ref{existAM}] We claim that a 
subsequence of the
$\Sigma^k$'s of Proposition~\ref{AMinpairs} satisfies the
requirements of Proposition~\ref{existAM}. 
Indeed fix $k\in \N$ and $r$ such that $\Inj
(M)>4r>0$. 
Since $(B_r (x), M\setminus B_{4r}(x))\in \co$
we then know that $\Sigma^N$ is $1/k$--a.m.\ in $M\setminus 
B_{4r} (x)$. Thus we have that 
\begin{eqnarray}
&&\mbox{either $\Sigma^k$ is $1/k$--a.m.\ on $B_r (y)$
for every $y$}\label{firstalt}\\
&&\mbox{or there is $x^k_r\in M$ s.t. $\Sigma^k$ is
$1/k$--a.m.\ on $M\setminus B_{4r} (x^k_r)$.}\label{secondalt}
\end{eqnarray}
If for some $r>0$ there 
exists a subsequence $\{\Sigma^{k(n)}\}$
satisfying (\ref{firstalt}), then we are done. Otherwise
we may assume that there are two sequences of
natural numbers $n\uparrow \infty, j\uparrow \infty$ and points $x^n_j$ such that
\begin{itemize}
\item For every $j$, and for $n$ large enough,
  $\Sigma^n$ is $1/n$--a.m.\ in $M\setminus B_{1/j} (x^n_j)$.
\item $x^n_j\to x_j$ for $n\uparrow \infty$ and 
$x_j\to x$ for $j\uparrow \infty$.
\end{itemize}
Thus for every $j$, the sequence $\{\Sigma^n\}$ is a.m.\
in $M\setminus B_{2/j} (x)$.
Of course if $U\subset V$ and $N$ is $\eps$--a.m.\ in $V$, then $N$
is $\eps$--a.m.\ in $U$. This proves that there exists a subsequence 
$\{\Sigma^j\}$ which
satisfies conditions  
\eqref{primai} and \eqref{primaii} for some positive function 
$r:M \to \rn{+}$.

It remains to show that an appropriate further subsequence satisfies 
\eqref{primaiii}. 
Each $\Sigma^j$ is smooth except at finitely many points. We denote by 
$P_j$ the set of 
singular points of $\Sigma^j$. After extracting another subsequence we 
can assume that 
$P_j$ is converging, in the Hausdorff topology, to a finite set $P$. If 
$x\in P$ and 
$\an$ is any annulus centered at $x$, then  $P_j \cap \an=\emptyset$ for 
$j$ large enough. 
If $x\not \in P$ and $\an$ is any (small) annulus centered at $x$ with outer 
radius less than 
$\d (x, P)$, then $P_j \cap \an=\emptyset$ for $j$ large enough. Thus, 
after possibly 
modifying the function $r$ above,
the sequence $\{\Sigma^j\}$ satisfies \eqref{primai}, \eqref{primaiii},  
and \eqref{primaii}.
\end{proof}

\section{Regularity for the replacements}\label{s:reg1}

We will now define a notion of a ``good replacement'' for stationary
varifolds and prove that the existence of (sufficiently many) 
replacements for a stationary varifold implies that it is a smooth 
minimal surface; see
Proposition~\ref{reg}. In 
section 6 we will show that the varifold $V$ of
Proposition~\ref{existAM} satisfies the hypotheses of
Proposition~\ref{reg} and thus is smooth.

\begin{definition}\label{repldef} Let $V\in \va$ be stationary and
$U\subset M$ be an open subset.
A stationary varifold $V'\in \va$ is said to be 
a {\em replacement for $V$ in $U$}\/ if \eqref{secondai} and
\eqref{secondaii} below hold.
\begin{eqnarray}
&&\mbox{$V'=V$ on $G(M\setminus \ov{U})$ and $\|V'\| (M)=\|V\| (M)$.}
\label{secondai}\\
&&\mbox{$V\res U$ is a stable minimal surface $\Sigma$
with  $\overline{\Sigma}\setminus \Sigma\subset \partial U$.}\label{secondaii}
\end{eqnarray}
\end{definition}

\begin{definition}\label{goodrepl} Let $V$ be a stationary
  varifold and $U\subset M$ be an open subset. 
We say that $V$ has the {\em good replacement property}\/ in $U$ if
(a), (b), and (c) below hold.
\begin{itemize}
\item[(a)] There is a positive function $r:U\to \rn{}$ such that for
  every annulus $\an\in \An_{r(x)} (x)$ there is a replacement for $V'$
  in $\an$.
\item[(b)] The replacement $V'$ has a replacement $V''$ in any
  $\an\in \An_{r(x)} (x)$ and in any $\an\in \An_{r'(y)} (y)$ 
(where $r'$ is positive).
\item[(c)] $V''$ has a replacement $V'''$ in any
  $\an\in \An_{r''(y)} (y)$ (where $r''>0$).
\end{itemize}
If $V$ and $V'$ are as above, then we
will say that $V'$ is a {\em good replacement}\/ and $V''$ a {\em
good further replacement}.  
\end{definition}

This section is devoted to prove the following:

\begin{propos}\label{reg} Let $G$ be an open subset of $M$. If
$V$ has the good replacement property in $G$, then $V$ is a
(smooth) minimal surface in $G$.
\end{propos}

In the proof Proposition~\ref{reg} we need the two 
technical Lemmas~\ref{PittsPr} and~\ref{l:tec}, 
stated and proved in Appendix~\ref{App2}. Note that
Lemma~\ref{PittsPr} is just a weak
version (in the framework of varifolds) of the classical maximum 
principle for minimal surfaces. 
As a first step towards the proof of Proposition~\ref{reg}
we have the following:

\begin{lemma}\label{rect} Let $U$ be an open subset of $M$ 
and $V$ a stationary
varifold in $U$.
If there exists a positive function $r$ on $M$ such that $V$ has 
a replacement in any annulus $\an\in \An_{r(x)} (x)$, then $V$
is integer rectifiable. Moreover, $\theta (x, V)\geq
1$ for any $x\in U$ and any tangent cone to $V$ in $x$ is an
integer multiple of a plane.
\end{lemma}

\begin{proof} Since $V$ is stationary, the monotonicity formula
  \eqref{e:MonFor} gives  a constant $C_M$ such that 
\begin{equation}\label{monot}
\frac{\|V\| (B_\sigma (x))}{\sigma^2}\;\leq\; C_M 
\frac{\|V\| (B_\rho (x))}{\rho^2} \qquad \forall \sigma<\rho<\Inj (M)\;\;
\mbox{and}\;\; \forall x\in M\,.
\end{equation}
Fix $x\in \supp (\|V\|)$ and $r<r(x)$ so that $4r$ is smaller than the
convexity radius. Replace 
$V$ with $V'$ in $\an (x, r, 2r)$. We claim that
$\|V'\|$ cannot be identically $0$ on $\An (x,r,2r)$. Assume it was;
since $x\in \supp (\|V'\|)$, there would be a $\rho\leq r$ such that
$V'$ ``touches'' $\partial B_\rho$ from the
interior. More precisely, there would exist $\rho$ and $\eps$
such that $\supp \|V'\|\cap \partial B_\rho (x)\neq \emptyset$ and
$\supp \|V'\|\cap \An (x, \rho, \rho+\eps)=\emptyset$. Since
$B_\rho (x)$ is convex this would
contradict Lemma \ref{PittsPr}. Thus $V'\res \an (x,r,2r)$ is
a non--empty smooth surface and so there is
$y\in \an(x, r, 2r)$ with $\theta (V', y)\geq 1$. Using
\eqref{monot} we get
\begin{equation}\label{estimate}
\frac{\|V\| (B_{4r} (x))}{16 r^2}\;=\;
\frac{\|V'\| (B_{4r} (x))}{16 r^2}\;\geq\;\frac{C_M \|V'\| (B_{2r} (y))}{16
  r^2}\;\stackrel{(\ref{monot})}{\geq}\; \frac{\pi C_M}{4}.
\end{equation}
Hence, $\theta (x, V)$ is bounded uniformly from below on $\supp
(\|V\|)$ and applying Theorem~\ref{Rectif} we conclude that
$V$ is rectifiable. 

We next prove that $V$ is {\em integer}\/ rectifiable. We use the
notation of Definition~\ref{cones}. Fix $x\in \supp
(\|V\|)$, a stationary
cone $C\in TV (x,V)$, and a sequence $\rho_n\downarrow 0$ such that
$V^x_{\rho_n} \to C$. Replace $V$ by $V'_n$ in
$\an (x, \rho_n/4, 3\rho_n/4)$ and set $W'_n=
(T^x_{\rho_n})_\sharp V'_n$. After possibly passing
to a subsequence, we can
assume that $W'_n \to C'$, where $C'$ is a stationary varifold. The following
properties of $C'$ are trivial consequences of the definition of
replacements:
\begin{eqnarray} 
&&\mbox{$C'=C$ in $\B_{1/4}(x)\cup \an (x, 3/4, 1)$,}\label{euno}\\ 
&&\mbox{$\|C'\| (\B_\rho)= \|C\| (\B_\rho)$ if $\rho\in ]0,1/4[\cup
  ]3/4,1[$}.\label{edue}
\end{eqnarray}
Since $C$ is a cone, in view of (\ref{edue}) we have
\begin{equation}\label{cone}
\frac{\|C'\| (\B_\sigma)}{\sigma^2}\;=\; \frac{\|C'\|
    (\B_\rho)}{\rho^2}\qquad \forall \sigma, \rho\in ]0,1/4[\cup ]3/4,1[.
\end{equation}
Hence, the stationarity of $C'$ 
and the monotonicity formula imply that $C'$ is a cone.
By \eqref{e:curvaclaim}, $W'_n$
converge to a stable embedded minimal surface in $\an (x, 1/4, 3/4)$.
This means that $C'\res \an(x,1/4, 3/4)$ is an
embedded minimal cone in the classical sense and hence it is supported
on a disk containing the origin. This forces $C'$ and $C$ to coincide
and be an integer multiple of the same plane.
\end{proof}

\begin{proof}[Proof of Proposition~\ref{reg}] The strategy of the
proof is as follows.
Fix $x\in M$, a good replacement $V'$ for $V$ in 
$\an(x, \rho, 2\rho)$, and
let $\Sigma'$ be 
the stable minimal surface given by $V'$ in $\an(x, \rho,
2\rho)$. Consider $t\in ]\rho, 2\rho[$, $s<\rho$
and the replacement $V''$ of $V'$ in $\an(x, s, t)$, which in this 
annulus coincides with a smooth minimal surface $\Sigma''$. 
In step 2 we will prove that, for $\rho$ sufficiently small
and for an appropriate choice of $t$, then $\Sigma''\cup \Sigma'$ is a 
smooth surface.
Letting $s\downarrow 0$ we get a minimal surface 
$\Sigma\subset B_\rho (x)\setminus \{x\}$ 
such that every $\Sigma''$ constructed as above is a subset of $\Sigma$. 
Loosely speaking, any replacement of $V'$ will coincide with $\Sigma$
in the annular region where it is smooth. 

Now, fix a $z$ which belongs to $\supp (\|V\|)$ and such that $V$ intersects
$\partial B_s (x)$ ``transversally'' in $z$.
If we consider a replacement $V''$ of $V'$ in $\an (x, s, \rho)$, then
$z$ will belong to the closure of the minimal surface $\Sigma''=V''\res \an
(x, s, \rho)$. The discussion above gives that
$z\in \Sigma$. Lemma~\ref{l:tec} implies 
that  ``transversality'' to the spheres centered at $x$ in a dense subset of 
$(\supp (\|V\|))\cap B_\rho (x)$.
Thus in step 3 we conclude that 
$$
(\supp (\|V\|))\cap B_\rho (x)\setminus \{x\}\;\subset\; \Sigma.
$$ 
Since 
$\haus^2 (\Sigma\cap B_\rho (x))=\|V\| (B_\rho (x))$, 
then $V=\Sigma$ in $B_\rho (x)$. Step 4 concludes the proof by showing  
that $x$ is a removable singularity for $\Sigma$.

The key fact that $\Sigma''$ and $\Sigma'$ can be ``glued'' smoothly together
is a consequence of 
the curvature estimates for stable minimal surfaces combined with
the characterization of the tangent cones given in Lemma~\ref{rect}.
These two ingredients will be used to prove that $\Sigma''$ is (locally) a 
Lipschitz graph nearby $\partial B_t (x)$; 
thus allowing us to apply standard theory of Elliptic PDE.

\medskip

\noindent
{\bf Step 1: The set up.}

Fix $x$, $V$, $V'$, $V''$, $\Sigma'$, $\Sigma''$, $\rho$, $s$, and $t$ 
as above.
We require that $2\rho$ is less than the convexity radius of $M$ and 
that $\Sigma'$ intersects 
$\partial B_t (x)$ transversally. Fix a point 
$y\in \Sigma'\cap \partial B_t (x)$ and
a sufficiently small radius $r$, so that $\Sigma'\cap B_r (y)$ 
is a disk and $\gamma=
\Sigma'\cap \partial B_t (x)\cap B_r (y)$ is a smooth arc.

Let $\zeta: B_r (y)\to \B_1$ be a diffeomorphism such that 
$$
\zeta (\partial B_t (x))\subset \{z_1=0\} \qquad \mbox{and} \qquad
\zeta(\Sigma'')\subset \{z_1>0\}\, ,
$$
where $z_1,z_2,z_3$ are orthonormal coordinates on $\B_1$;
see Fig.~\ref{f:2}. We will also assume that 
$\zeta (\gamma)=\{(0,z_2, g'(0,z_2))\}$ and
$\zeta (\Sigma')\cap \{z_1\leq 0\}=\{(z_1, z_2, g' (z_1, z_2))\}$ 
where $g'$ is smooth.

\begin{figure}[htbp]
\begin{center}
    \input{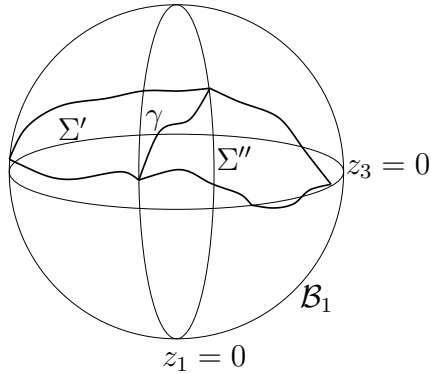}
    \caption{The surfaces $\Sigma'$ and $\Sigma''$ and the curve
      $\gamma$ in $\B_1$.} 
    \label{f:2}
\end{center}
\end{figure}

Note the following elementary facts:
\begin{itemize}
\item Any kind of estimates (like curvature or area bounds or monotonicity) 
for a minimal surface $\Sigma\subset B_r(y)$
translates into similar estimates for the surface $\zeta (\Sigma)$.
\item Varifolds in $B_r (y)$ are push--forwarded to varifolds in $\B_1$ and 
there is a natural correspondence 
between tangent cones to $V$ in $\xi$ and tangent cones to $\zeta_\sharp V$ in 
$\zeta (\xi)$.
\end{itemize}
By slight abuse of notation, we use the same symbols (e.g. $\gamma$, $V'$, 
$\Sigma'$) for both the objects of $B_r (y)$ and
their images under $\zeta$.

\medskip

{\bf Step 2: Graphicality; gluing $\Sigma'$ and $\Sigma''$ smoothly together. }

The varifold $V''$ consists of $\Sigma''\cup \Sigma'$ in $B_r
(y)$. Moreover, 
Lemma~\ref{rect} applied to $V''$ gives that $TV(z, V'')$ is a family 
of (multiples of) $2$--planes.
Fix $\overline{z}\in\gamma$. Since $\Sigma'$ is regular and transversal to 
$\{z_1=0\}$ in $\overline{z}$, each plane $P\in TV(\overline{z},V'')$ 
coincides with the half plane $T_{\overline{z}} \Sigma'$ in $\{z_1<0\}$. Hence 
$TV (\overline{z}, V'')=\{T_{\overline{z}} \Sigma'\}$.
Let $\tau (\overline{z})$ be the unit normal to the graph of $g'$
$$
\tau (\overline{z})\;=\; \frac{(-\partial_1 g' (0, \overline{z}_2), 
-\partial_2 g' (0, \overline{z}_2), 1)}
{\sqrt{1+|\nabla g'0,\overline{z}_2)|^2}}
$$
and let $R^{\overline{z}}_r: \rn{3}\to \rn{3}$ be the 
dilatation of $3$-space defined by
$$
R^{\overline{z}}_r (z) \;=\; \frac{z-\overline{z}}{r}\,.
$$
Since $TV(\overline{z}, V'')= \{T_{\overline{z}} \Sigma'\}$, the surfaces 
$\Sigma_r= R^{\overline{z}}_r (\Sigma'')$ converge to the 
half plane $HP= \{\tau (\overline{z})\cdot v =0, v_1>0\}$ --- half of the 
plane $\{\tau (\overline{z})\cdot v =0\}$. This
convergence implies that
\begin{equation}\label{e:L1}
\lim_{z\to \overline{z}, z\in \Sigma''} 
\frac{|(z-\overline{z})\cdot \tau (\overline{z})|}{|\overline{z}-z|}=0\,.
\end{equation}
Indeed assume that \eqref{e:L1} fails; 
then there is a sequence $\{z_n\}\subset \Sigma''$ such that 
$z_n \to \overline{z}$ and
$|(z_n-\overline{z})\cdot \tau (\overline{z})|\geq k |z_n-\overline{z}|$ 
for some $k>0$. Set $r_n=|z_n -\overline{z}|$. There exists a constant 
$k_2$ such that $\B_{2k_2 r_n} (z_n)\cap HP=\emptyset$. Thus 
$\dist (HP,\B_{k_2 r_n} (z_n))\geq k_2 r_n$.
Since $\Sigma''$ is regular in $z_n$ we get by the monotonicity formula that
$$
\|V''\| (\B_{k_2 r_n} (z_n))\geq C k^2_2 r_n^2 \qquad \mbox{where $C$ 
depends on $\zeta$.}
$$
This contradicts the fact that $HP$ is the only element of 
$TV(\overline{z}, V'')$.
Note also that the convergence of \eqref{e:L1} is uniform for 
$\overline{z}$ in compact subsets of $\gamma$.
The argument is explained in Fig.~\ref{f:6}.

\begin{figure}[htbp]
\begin{center}
    \input{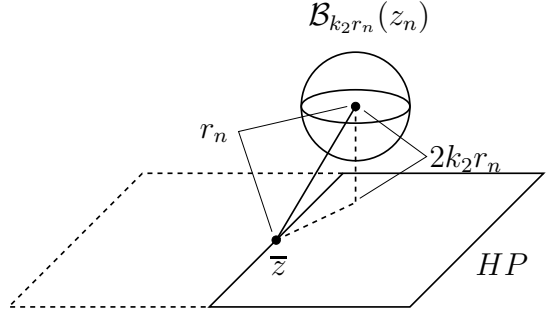}
    \caption{If $z_n\in \Sigma''$ is far from the plane $HP$, 
the monotonicity formula gives a ``good amount'' of the varifold $V''$ which
    lives far from $HP$.}
    \label{f:6}
\end{center}
\end{figure}

Let $\nu$ denote the smooth unit vector field to $\Sigma''$ such that 
$\nu\cdot (0,0,1)\geq 0$.
We next use the stability of $\Sigma''$ to show that 
\begin{equation}\label{e:L2}
\lim_{z\to \overline{z}, z\in \Sigma''} \nu (z) \;=\; \tau
(\overline{z})\, .
\end{equation}
Indeed let $\sigma$ be the plane $\{(0,\alpha, \beta), \, \alpha,
\beta\in \rn{}\}$, 
assume that $z_n\to \overline{z}$ and set $r_n=\dist (z_n,\sigma)$. Define
the rescaled surfaces
$\Sigma^n = R^{z_n}_{r_n} (\Sigma''_n \cap \B_{r_n} (z_n))$. Each 
$\Sigma^n$ is a stable 
minimal surface in $\B_1$, and hence, after possibly passing to a 
subsequence,  
$\Sigma^n$ converges smoothly in 
$\B_{1/2}$ to a minimal
surface $\Sigma^\infty$ (by \eqref{e:curvaclaim}). By
\eqref{e:L1}, we have that $\Sigma^\infty$ is the disk 
$T_{\overline{z}}\Sigma'\cap \B_{1/2}$. 
Thus the normals to $\Sigma^n$ in $0$, which are given 
by $\nu (z_n)$, converge to $\tau (\overline{z})$; see Fig.~\ref{f:7}. 
It is easy 
to see that the convergence in \eqref{e:L2} is uniform on compact subsets 
of $\gamma$.

\begin{figure}[htbp]
\begin{center}
    \input{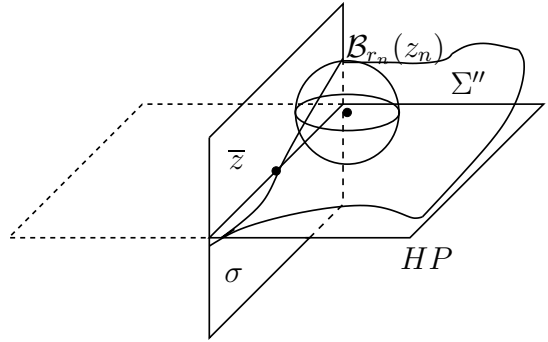}
    \caption{If we rescale $\B_{r_n} (z_n)$, then we find a sequence of 
stable minimal surfaces $\Sigma^n$ which converge to the half--plane $HP$.}
    \label{f:7}
\end{center}
\end{figure}

Hence, for each $\overline{z}\in \gamma$, there exists $r>0$ and a function 
$g''\in C^1(\{z_1\geq 0\})$ such that
\begin{eqnarray*}
&&\Sigma''\cap B_r (\overline{z})\;=\;\{(z_1, z_2, g''(z_1, z_2)),
z_1>0\}\, ,\\ 
&&g''(0, z_2)=g'(0, z_2)\, , \qquad \mbox{and} 
\qquad \nabla g'' (0, z_2)= \nabla g' (0, z_2)\,.
\end{eqnarray*} 
In the coordinates $z_1, z_2, z_3$, the minimal surface equation 
yields a second order uniformly 
elliptic equation for $g'$ and
$g''$. Thus the classical theory of elliptic PDE 
gives that $g'$ and $g''$ are restrictions of a unique smooth function $g$.

\medskip

\noindent
{\bf Step 3: Regularity of $V$ in the punctured ball.}

Let $\Sigma'$ and $\Sigma''$ be as in the previous step.
We will now show that :
\begin{equation}\label{e:claims3}
\mbox{If $\Gamma$ is a connected component of $\Sigma''$, then 
$\overline{\Gamma} \cap \Sigma' \cap \partial B_t (y)\neq \emptyset$.}
\end{equation}
Indeed assume that for some $\Gamma$ equation \eqr{e:claims3} fails.  
Since $t$ is assumed to be 
less than the convexity radius we have by the maximum 
principle that 
$\overline{\Gamma}\cap \partial B_t (x) \neq \emptyset$. Fix $z$ in
$\overline{\Gamma}\cap \partial B_t (x)$. 
If \eqref{e:claims3} were false, then the
varifold $V''$ would ``touch'' $\partial B_t (x)$ in $z$ from the
interior. More precisely, there would be an $r>0$ such that
$$
z\in \supp (\|V''\|)\qquad \mbox{and}\qquad  
(B_r (y)\cap \supp (\|V''\|)) \subset \overline{B}_t (x)\,.
$$

This contradicts Lemma~\ref{PittsPr}; 
thus \eqref{e:claims3} holds.

Let $t, \rho$ be as in the first paragraph of Step 1. Step 2 and 
\eqref{e:claims3} imply the following:
\begin{eqnarray}
&&\mbox{if $s<\rho$, then $\Sigma'$ can be extended to}\nonumber\\
&&\mbox{a surface $\Sigma_s$ in $\an(x, s, 2\rho)$}\qquad\\ 
&&\mbox{if $s_1<s_2< \rho$, then $\Sigma_{s_1}=\Sigma_{s_2}$ in 
$\an (x, s_2, 2\rho)$.} 
\end{eqnarray}
Thus $\Sigma=\bigcup_s \Sigma_s$ is a stable minimal surface 
$\Sigma$ with $\overline{\Sigma}\setminus \Sigma\subset 
(\partial B_{2\rho} (x)\cup\{x\})$, i.e. $\Sigma'$ can be continued up
to $x$ (which, in principle, could be a singular point).
  
We will next show that $V$ coincides with $\Sigma$ in $B_\rho (x)\setminus
\{x\}$. Recall that $V=V'$ in $B_\rho (x)$. Fix 
$$
y\in (\supp (\|V\|))\cap B_\rho (x) \setminus \{x\} \qquad 
\mbox{and set $s= d (y,x)$}.
$$
We first prove that if $TV(y, V)$ consists of a (multiple of a)
plane $\pi$ transversal 
to $\partial B_s (x)$, then $y$ belongs to $\Sigma$.
Consider the replacement $V''$ of $V'$ in $\an (x, s, t)$ and split $V''$ 
into the three varifolds 
$$
\begin{array}{lllll}
V_1 &=& V''\res B_s (x)&=&V\res B_s (x)\, ,\\
V_2 &=& V''\res \an (x, s, 2\rho)&=&\Sigma \cap \an (x, s, 2\rho)\, ,\\
V_3 &=&V''- V_1 -V_2\, .&&\\
\end{array}
$$
By Lemma~\ref{rect}, the set $TV(y, V'')$ 
consists of planes and since $V_1=V\res B_s (x)$, all these planes
have to be multiples of $\pi$. Thus 
$y$ is in the closure of  $(\supp (\|V''\|))\setminus \overline{B}_s
(x)$, which implies $y\in \overline{\Sigma}_t\subset \Sigma$. 

Let $T$ be the set of points $y\in B_\rho (x)$ such that
$TV(y, V)$ consists of a (multiple of) a plane transversal to 
$\partial B_{d(y,x)} (x)$. 
Lemma~\ref{l:tec} gives that $T$ is dense in $\supp (\|V\|)$.
Thus 
$$
(\supp (\|V\|))\cap B_\rho (x)\setminus \{x\}\;\subset\; \Sigma\, .
$$ 
Property \eqref{secondai} of replacements implies
$\haus^2 (\Sigma\cap B_\rho (x))= \|V\| (B_\rho (x))$. Hence $V=\Sigma$ 
on $B_\rho (x)\setminus \{x\}$.

\medskip

\noindent
{\bf Step 4: Regularity in $x$.}

We will next show that $\Sigma$ is smooth also in $x$, i.e.\ that
$x$ is a removable singularity for $\Sigma$.
If $x\not \in \supp (\|V\|)$, then we are done. So assume that 
$x\in \supp (\|V\|)$. In the following we will use that, by
Lemma~\ref{rect}, every $C\in TV (x, V)$ is a multiple of a
plane. 

Map $B_t (x)$ into $\B_t (0)$ by the exponential map, use the
notation of Step 1, and set
$\Sigma_r = R^x_r (\Sigma)$. Every convergent subsequence 
$\{\Sigma_{r_n}\}$ converges to a plane in the sense of varifolds. 
The curvature estimates for stable minimal surfaces 
(see \eqref{e:curvaclaim})
gives that this convergence is actually smooth in $\B_1\setminus \B_{1/2}$. 
Thus, for $r$ sufficiently small, there exist natural numbers $N
(\rho)$ and $m_i (\rho)$ such that
$$
\Sigma\cap \an (x, \rho/2, \rho) \;=\;
\bigcup_{i=1}^{N (\rho)} m_i (\rho) \Sigma^i_\rho\, ,
$$ 
where each $\Sigma_\rho^i$ is
a Lipschitz graph over a (planar) annulus. Note also that
the Lipschitz constants are uniformly bounded, independently of $\rho$.

By continuity, the numbers $N(r)$ and $m_i (r)$ do not
depend on $r$. Moreover, if $s\in ]\rho/2, \rho[$, then each
$\Sigma_\rho^i$ can be continued through $\an (s/2, \rho/2, x)$ by a
$\Sigma_s^j$. Repeating this argument a countable number of times, 
we get $N$ minimal punctured disks $\Sigma^i$ with
$$
\Sigma\cap B_\rho (x)\setminus \{x\}\;=\;\bigcup_{i=1}^N m_i
\Sigma^i\, .
$$ 
Note that $x$ is a removable singularity for each
$\Sigma^i$. Indeed, $\Sigma^i$ is a stationary varifold in $B_\rho
(x)$ and $TV(x, \Sigma^i)$ consists of planes with multiplicity one.
This means that 
$$
\lim_{r\downarrow 0} \frac{\|V\| (B_r (x))}{\pi r^2} \;=\;1\, .
$$
Hence we can apply Allard's regularity theorem (see section 8 of \cite{All})
to conclude that $\Sigma^i$ is a graph in a sufficiently small ball
around $x$. Standard elliptic PDE theory gives that $x$ is a removable
singularity.

Finally, the maximum principle for minimal surfaces implies
that $N$ must be $1$. This completes the proof.
\end{proof}

\begin{remark} In the case at hand, there are other ways of
  proving that $x$ is a removable singularity. 
For example one could use the
existence of a conformal parameterization $u:
\CC\setminus \{0\}\to \Sigma^i\cap B_\rho (x)\setminus \{x\}$. The
minimality of $\Sigma^i$ gives that $u$ is an harmonic map. Since the
energy of $u$ is finite, we
can use theorem 3.6 in \cite{SU} to conclude that $u$ is smooth in $0$.
\end{remark}

\section{Construction of the replacements}\label{s:reg2}

In this section we conclude the proof of Theorem \ref{t:SS} by showing that
the varifold $V$ of Proposition \ref{existAM} is a smooth minimal
surface. 

\begin{theorem}\label{goal}
Let $\{\Sigma^j\}$ be a sequence of compact surfaces in
$M$ which converge to a stationary varifold $V$. If 
there exists a function $r:M\to \rn{+}$ such that 
\begin{itemize}
\item in every annulus of $\An_{r(x)} (x)$ and 
for $j$ large enough $\Sigma^j$ is a $1/j$--a.m. smooth surface in $\an$,
\end{itemize} 
then $V$ is a smooth minimal surface.  
\end{theorem}

To prove this theorem, we will show that 
$V$ satisfies the requirements of Proposition~\ref{reg}. Thus we 
need to construct good replacements for $V$, using the strategy outlined in 
Section~\ref{s:ov}.
In subsection~\ref{ss:MSY} we fix some notation and recall a
theorem of Meeks--Simon--Yau. 
In subsection~\ref{ss:con} we show how to construct the varifolds
$V^*$ which are our candidates for replacements. 
Subsections~\ref{ss:3} and~\ref{ss:4} prove the regularity of the 
$V^*$'s constructed in subsection~\ref{ss:con}. Finally, in 
subsection~\ref{ss:goal} 
we prove the last details needed to show that $V$ meets the 
requirements of Proposition~\ref{reg}.

\subsection{The result of Meeks--Simon--Yau}\label{ss:MSY}

\begin{definition}\label{minprinciple}
Let $\mathcal{I}$ 
be a class of isotopies of $M$ and $\Sigma\subset M$ a smooth embedded
surface. If $\{\varphi^k\}\subset \mathcal{I}$ and 
$$
\lim_{k\to\infty} 
\haus^2 (\varphi^k (1,\Sigma))\;=\;\inf_{\psi\in \mathcal{I}} 
\haus^2 (\psi (1,\Sigma))\, ,
$$
then we say that $\varphi^k(1,\Sigma)$ is a {\em minimizing
  sequence for Problem $(\Sigma, \mathcal{I})$}.
\end{definition}

We will need the following result from\cite{MSY}:

\begin{theorem}\label{theoMSY}{\rm [Meeks--Simon--Yau \cite{MSY}]} If
  $\{\Sigma^k\}$ is a minimizing sequence for Problem $(\Sigma, \Is
  (U))$ which converges to a varifold $V$,
  then there exists a stable minimal surface $\Gamma$ with 
  $\overline{\Gamma}\setminus \Gamma \subset \partial U$ and 
  $V=\Gamma$ in $U$. 
\end{theorem}

In \cite{MSY} Theorem~\ref{theoMSY} is proved for $U=M$.
  However, the theory developed there is local and
  can be extended in a straightforward way to cover the case at hand.

\subsection{Construction of replacements}\label{ss:con}
Let $V$ be as in Theorem~\ref{goal} and fix an annulus $\an\in
\An_{r(x)} (x)$. Set
\begin{eqnarray*}
\lefteqn{\Is_j (\an)\;=}&&\\
&&\bigl\{\psi\in \Is (\an)\bigr|\, \, 
\haus^2 (\psi_j (\tau, \Sigma^j))
\leq \haus^2 (\Sigma_j)+1/(8j)\quad \forall 
\tau\in [0,1]\bigr\}\, ,\\
\lefteqn{m_j=\inf_{\psi\in \Is_j} \haus^2 (\psi (1, \Sigma^j))\, .}&&\\
\end{eqnarray*}

Fix $j$. The following lemma implies that we can deform $\Sigma^j$ into
a sequence $\Sigma^{j,k}$ which is minimizing for Problem
$(\Sigma^j, \Is_j (\an))$ and converges, as $k\to \infty$,  
to a stable minimal surface in $\an$.

\begin{lemma}\label{stablemin}
If a sequence $\{\Sigma^{j,k}\}^k$ is minimizing for
Problem\\ 
$(\Sigma^j, \Is_j (\an))$ and converges
to a varifold $V^j$, then $V^j$ is a stable minimal surface in $\an$.
\end{lemma}

Lemma \ref{stablemin} will be proved in the next subsection. Here we use it for
constructing a replacement for $V$ in $\an$.

\begin{propos}\label{repl}
Let $V^j$ be the varifold of Lemma~\ref{stablemin}. Any $V^*$ which 
is the limit of a subsequence of
$\{V^j\}$ is a replacement for $V$.
\end{propos}

\begin{proof} Without loss of generality, we can assume that the 
  sequence $\{V^j\}$ converges to $V$.
Note that every $V^j$ coincides with $V$ in $M\setminus \ov{\an}$; thus the
same is true for $V^*$. Moreover, $\|V^j\| (M)\geq \haus^2
(\Sigma^j)-1/j$ since $\Sigma^j$ is a.m.
This gives that $\|V^*\| (M)=\|V\| (M)$. By
Lemma~\ref{stablemin} and 
\eqref{e:curvaclaim} we have that $V^*$ is a stable minimal 
surface in $\an$. 

To complete the proof we need to show that $V^*$ is stationary.
Since $V=V^*$ in $M\setminus \ov{\an}$, then 
$V^*$ is stationary in this open set. Hence it
suffices to prove that $V^*$ is stationary in an open annulus $\an'\in \An_r$
containing $\ov{\an}$. Choose such an $\an'$ and
suppose that $V^*$ 
is not stationary in $\an'$; we will show that this contradict that
$\{\Sigma^j\}$ is a.m.\ in $\an'$. Namely, suppose that for some vector
field $\chi$ supported in $\an'$ we have $\delta V^* (\chi)\leq
-C<0$. Let $\psi$ be the isotopy given by that ${\textstyle 
\frac{\partial \psi (t, x)}{\partial t}}=\chi (\psi (t,x))$ and set
$$
\begin{array}{lll}
V^* (t) &=& \psi(t)_\sharp V^*\, ,\\
V^j (t) &=& \psi(t)_\sharp V^j\, ,\\ 
\Sigma^{j,k} (t)&=& \psi (t, \Sigma^{j,k})\, .
\end{array}
$$
For $\eps$ sufficiently small, we have that
$$
[\delta V^* (t)] (\chi)\;\leq\; -\frac{C}{2} \qquad \mbox{for all $t<\eps$}.
$$
Since $V^j (t) \to V^* (t)$, there exists $J$ such that 
$$
[\delta V^j (t)] (\chi)\leq -\frac{C}{4} 
\qquad \mbox{for every $j>J$ and every $t<\eps$}.
$$
Moreover, since $\Sigma^{j,k} (t)\to V^j (t)$, 
for each $j>J$ there exists $K(j)$ with 
\begin{equation}\label{gaindiff}
[\delta \Sigma^{j,k} (t)] (\chi)\;\leq\; -\frac{C}{8}\qquad
\mbox{for all $t<\eps$ and all $k>K(j)$.}
\end{equation}
Integrating both sides of (\ref{gaindiff}) we get
\begin{equation}\label{gain}
\haus^2 (\psi (t, \Sigma^{j,k}))- \haus^2 (\Sigma^{j,k})\;\leq\; -\frac{t
  C}{8}
\end{equation}
Choose $j$ and $k$ 
sufficiently large so that $\eps C/8> 1/j$ and \eqref{gain} holds.
Each $\Sigma^{j,k}$ is isotopic to $\Sigma^j$ via an isotopy 
$\varphi^{j,k}\in \Is_j (\an)$.
By gluing $\varphi^{j,k}$ and $\psi$ smoothly together, we find a 
smooth isotopy $\Phi:[0, 1+\eps]\times M\to
M$ supported on $\an'$. In view of \eqref{gain}, $\Phi$ satisfies
\begin{eqnarray*}
\haus^2 (\Phi (t, \Sigma^j)) &\leq& \haus^2 (\Sigma^j)+1/(8j) \qquad
\forall t\in [0, 1+\eps]\, ,\\
\haus^2 (\Phi (1+\eps, \Sigma^j) &<& \haus^2 (\Sigma^j)-1/j\, ,
\end{eqnarray*}
which give the desired contradiction and prove the proposition.
\end{proof}

\subsection{Proof of Lemma~\ref{stablemin}}\label{ss:3}
Without loss of generality we may assume that $j=1$ and use
$V'$, $\Sigma$ and $\Sigma^k$ in place of $V^j$, $\Sigma^{j,k}$ and
$\Sigma^j$. Clearly $V'$ is stationary and stable in $\an$, by
its minimizing property. 
Thus we need only prove that $V'$ is regular. The proof of this uses 
Theorem~\ref{theoMSY} and the following:

\begin{lemma}\label{applyMSY}
Let $x\in \an$ and assume that
$\{\Sigma^k\}$ is minimizing for Problem $(\Sigma, \Is_1 (\an))$. 
There exists $\eps>0$ such that, for $k$ sufficiently
large, the following holds: 
\begin{itemize}
\item[(Cl)] For any $\varphi\in \Is (B_\eps (x))$ with
$\haus^2 (\varphi (1, \Sigma^k))\leq \haus^2 (\Sigma^k)$, 
there exists an isotopy $\Phi\in \Is (B_\eps)$ such that
\begin{eqnarray}
&&\Phi(1, \cdot)\;=\; \varphi (1, \cdot)\,,\label{itcoincides}\\
&&\haus^2 (\Phi (t, \Sigma^k))\;\leq\; \haus^2 (\Sigma^k)+1/8\label{itisOK}.
\end{eqnarray}
\end{itemize}
Moreover, $\eps>0$ can be chosen so that (Cl) holds for
{\em any} sequence $\{\tilde{\Sigma}^k\}$ which is minimizing 
for Problem $(\Sigma, \Is_1 (\an))$
and with $\Sigma^j=\tilde{\Sigma}^j$ 
on $M\setminus \ov{B}_\eps (x)$.
\end{lemma}

Lemma~\ref{applyMSY} will be proved in the next subsection. We now return
to the proof of Lemma~\ref{stablemin}. We will use
Proposition~\ref{reg}. Hence, once again we need to construct
replacements for a varifold, which this time is $V'$. 
We divide the proof into two steps. The first one is the basic
construction of replacements for $V'$. The second shows that the
replacements satisfy (b) and (c) in 
Definition~\ref{goodrepl}

\medskip
 
{\bf Step 1} Fix $x\in \an$ and $\eps>0$ such that Lemma~\ref{applyMSY}
holds. Fix any annulus 
$$
\an^*\;=\;\an(x, \tau, t)\;\subset\; B_\eps (x)\;\subset\; \an
$$ 
and consider a minimizing sequence $\{\ov{\Sigma}^{k,l}\}^l$ for
Problem $(\Sigma^k, \Is (\an^*))$. Lemma~\ref{applyMSY} 
implies that, for $k$ sufficiently large,
$\ov{\Sigma}^{k,l}$ can be constructed from $\Sigma$ via an isotopy
of $\Is_1 (\an)$.
Thus if we let $W^k$ be
the varifold limit of $\ov{\Sigma}^{k,l}$ and $W$ the limit of
$W^k$, then we have that $\|W\| (M)=\|V'\| (M)$. 

Let $\{\ov{\Sigma}^{k,l(k)}\}$ be a subsequence which converges
to $W$. By the discussion above the subsequence is a minimizing 
sequence for Problem $(\Sigma, \Is_1 (\an))$. 
Hence $W$ is stationary in
$\an$. Moreover, by Theorem~\ref{theoMSY}, every $W^k$ is a stable
minimal surface in the annulus $\an^*$: the 
curvature estimates (see \eqref{e:curvaclaim})
give that $W$ is a stable minimal surface in $\an^*$.
Hence $W$ is a replacement for $V'$.

\medskip

{\bf Step 2} Summarizing we have proven in Step 1 that for any $y\in \an$ 
there exists $r(y)>0$ such that
in the class of annuli $\An_{r(y)} (y)$ we can construct replacements. 
In order to complete the proof we have to check that the replacements
so constructed
satisfy all the technical requirements of Proposition~\ref{reg}. Thus 
define $W$ as in Step 1.
Since $\{\ov{\Sigma}^{k,l(k)}\}$ is a
minimizing sequence for Problem $(\Sigma, \Is_1 (\an))$, 
in all the arguments of Step 1 we can use $W$ in place of $V'$. Thus 
for every $y\in \an$, $W$ has the replacement property for a class of
annuli centered at $y$. 
This shows the second part of (b) in Definition~\ref{goodrepl}

We still have to settle the first part of (b) in 
Definition~\ref{goodrepl}, i.e.\ 
that if the $W$ constructed in Step 1
replaces $V'$ in an annulus of $\An_{r(x)} (x)$,
then $W$ has the replacement property on the {\em whole}\/
collection of annuli $\An_{r(x)} (x)$. 
Note that our $r(x)=\eps$, where $\eps$ is given by Lemma~\ref{applyMSY}.
Every $\ov{\Sigma}^{k,l(k)}$
coincides with $\Sigma^k$ in $M\setminus B_\eps (x)$ and 
$\{\ov{\Sigma}^{k,l(k)}\}$
is minimizing for Problem $(\Sigma, \Is_1 (\an))$. Thus the last line of
Lemma~\ref{applyMSY} applies and again we can argue as in Step 1 with
$W$ in place of $V'$. We conclude that also $W$ has a replacement in every 
annulus of $\An_\eps (x)$.

Condition (c) in Definition~\ref{goodrepl} follows from similar arguments.
Summarizing, $V'$ and all its replacements just constructed
satisfy the requirements of Proposition~\ref{reg}. Hence $V'$ is a 
smooth surface in $\an$.
\qed 

\noindent
\subsection{Proof of Lemma~\ref{applyMSY}}{\bf Step 1: Small area slices.}
\label{ss:4}
Let $x$, $\eps$, $\an$, $\Sigma^k$ and $\varphi$ be as in the statement of the 
lemma.
We know that $\Sigma^k$ converges to a
varifold $V'$ which is stationary in $\an$ and 
has mass $m_0$. Thus the monotonicity formula gives that 
  $$
\|V'\| (\an (x, \eps, 2\eps))\;<\; 4 m_0 \eps^ 2\, .\qquad 
$$ 
Therefore, if $k$ is
sufficiently large we get that  
$\haus^2 (\Sigma^k\cap B_{2\eps} (\delta))< 5 m_0 \eps^2$. 
Applying the coarea formula we have that, for every
  such $k$, there exists an interval $I_k\subset ]\eps, 2\eps[$ such that 
\begin{equation}\label{small1}
\leb^1 (I_k)>0\quad \mbox{and} \quad
\haus^1 (\Sigma^k\cap \partial B_{\tau} (x))\;\leq\;
10m_0\eps \quad \mbox{for all $\tau\in I_k$}\, .
\end{equation}
Thus, applying Sard's Theorem to the
function $\d (\cdot, x)$ on $\Sigma^k$, we can
find $\tau_k\in ]\eps, 2\eps[$ such that
\begin{equation}\label{smallslice}
\haus^1 (\Sigma^k\cap \partial B_{\tau_k} (x))\;<\;
10 m_0\eps \quad \mbox{and }\quad \Sigma^k \text{ is transversal to
  $\partial B_{\tau_k} (x)$}\, .
\end{equation}
Moreover, the smoothness of $\Sigma^k$ implies that we can choose a small
interval $]\sigma_k, s_k[$ with $\eps<\sigma_k$ and so that 
\eqr{smallslice} holds for every $\tau_k\in ]\sigma_k, s_k[$.

\medskip

\noindent
{\bf Step 2: Radial deformations.}

For $\eta>0$ we denote by $O_\eta$ the usual radial deformation 
of Euclidean $3$--space given by
$O_\eta (x)= \eta x$. If both $r$ and $r\eta$ are less than $\Inj
(M)/2$, then we define the diffeomorphism
$$
I_\eta: B_r (x)\,\to\, B_{\eta r} (x)\qquad \mbox{by}\;\;\exp\circ
O_\eta \circ \exp^{-1}. 
$$ 
By the smoothness of $M$ there exists $\mu>0$ such that, 
for any surface $\Gamma\subset M$,
\begin{eqnarray}
\haus^2 (I_\eta (\Gamma) \cap B_{\eta r} (x))&\leq& \mu\, \eta^2\,
\haus^2 (\Gamma \cap B_r (x))\, ,\label{stima1}\\
\haus^1 (I_\eta (\Gamma)\cap \partial B_{\eta r} (x))&\leq& \mu\,\eta\,
\haus^1 (\Gamma\cap \partial B_r (x))\nonumber\\ 
&&\quad \mbox{if $\Gamma$ is transversal to 
$\partial B_r (x)$.}\label{stima2}\\
\haus^2 (\Gamma\cap B_r (x)) &\leq& 
\mu \int_0^r \haus^1 (\Gamma\cap \partial B_\rho
(x))\, d\rho.\label{coarea}
\end{eqnarray}
Fix $k$ and choose $\eps, \sigma_k$ and $s_k$ as in the
previous step. In the current step we use $I_\eta$ to construct a 
smooth $\psi:[0,1] \times M \to M$ such that 
\begin{itemize}
\item[-] For every $\delta>0$, $\psi|_{[0, 1-\delta]\times M}$ 
is a smooth isotopy supported in $\Is (B_{s_k} (x))$;
\item[-] $\psi|_{\{1\}\times M}$ ``squeezes'' the ball $B_{\sigma_k} (x)$ to
  the point $\{x\}$ and\\ 
``stretches'' the annulus $\an(x, \sigma_k, s_k)$ to
the ball $B_{s_k} (x)$; 
\item[-] For some constant $C$ depending on $\mu$ we have (see Fig.~\ref{f:4}) 
\begin{equation}\label{core}
\haus^2 (\psi (t, \Sigma^k))\;\leq\;\haus^2 (\Sigma^k) +
 C\eps^2\,.
\end{equation}
\end{itemize}
We construct $\psi$ explicitly. 
We first choose a nondecreasing 
smooth $f:[0,1]\times [0,1] \to [0,1]$ such that
$$
f(t,r)\;=\;
\left\{\begin{array}{ll}
(1-t) & \mbox{if $r\in [0, \sigma_k]$,}\\
1 & \mbox{if $r\in [s_k, 1]$,}
\end{array}\right.
$$
and then set
$$
\psi (t, y)\;=\;
\left\{\begin{array}{ll}
y & \mbox{if $\d (y,x)\geq s_k$}\, ,\\
I_{f(t, \d (x,y))} (y) & \mbox{otherwise.}
\end{array}\right.
$$

\begin{figure}[htbp]
\begin{center}
    \input{Pic4.pstex_t}
    \caption{The map $\psi (t, \cdot): M\to M$}
    \label{f:4}
\end{center}
\end{figure}

We will only prove \eqr{core}, since the other
properties are easy to check. First of all, since $\psi (t,y)=y$
on $M\setminus B_{\sigma_k} (x)$, we have
\begin{equation}\label{fuori}
\haus^2 (\psi (t, \Sigma^k) \cap (M\setminus B_{2\eps}
(x))\;=\;\haus^2 (\Sigma^k \cap (M\setminus B_{2\eps} (x)).
\end{equation}
By \eqr{stima1}
\begin{equation}\label{smallball}
\haus^2 (\psi (t, \Sigma^k)\cap B_{(1-t)\sigma_k} (x))\;\leq\;\mu (1-t)^2
\haus^2 (\Sigma^k\cap B_{\sigma_k})\; \leq \; 10 \mu m_0 \eps\, .
\end{equation}
To estimate the remaining portion of $\psi (t, \Sigma^k)$ we use
(\ref{coarea}) and get
\begin{eqnarray}
\lefteqn{\haus^2 (\psi (t, \Sigma^k)\cap \an (x, (1-t)\sigma_k, s_k))\;\leq}
&&\nonumber\\
&&\mu \int_{(1-t)\sigma_k}^{s_k} \haus^1 
(\psi (t, \Sigma^k)\cap \partial B_\rho (x)) d\rho.\label{almostOK}
\end{eqnarray}
Note that for $\rho\in ((1-t)\sigma_k, s_k)$
there exists $\tau\in (\sigma_k, s_k)$ and $\eta\leq 1$ such that 
$$
\psi (t, \Sigma^k)\cap \partial B_{\rho}\;=\;I_{\eta} (\Sigma^k\cap
\partial B_\tau (x))\, .
$$
Thus, by (\ref{stima2}) and (\ref{smallslice}) we have
\begin{equation}\label{stima3}
\haus^1 (\psi (t, \Sigma^k)\cap \partial B_{\rho})\;\leq\;
\mu \haus^1 (\Sigma^k\cap \partial B_\tau (x))\;\leq\; 10 \mu m_0 \eps\,.
\end{equation}
Inserting (\ref{stima3}) in (\ref{almostOK}) we find that
\begin{equation}\label{OK}
\haus^2 (\psi (t, \Sigma^k)\cap \an (x, (1-t)\sigma_k, s_k))\;\leq\;
20 \mu^2 m_0 \eps^2.
\end{equation} 
Equations (\ref{fuori}), (\ref{smallball}), and (\ref{OK}) yield
\begin{equation}\label{core2}
\haus^2 (\psi (t, \Sigma^k))\;\leq\;\haus^2 (\Sigma^k) + C\eps^2 \,.
\end{equation}

\medskip

\noindent
{\bf Step 3: The conclusion.}

Choose $\eps$ such that $\mu C\eps^2 <1/32$ 
and $K$ such that :
\begin{itemize}
\item We can construct the $\psi$ of the
previous step with (\ref{core2}) valid for any $k>K$;
\item $\haus^2 (\Sigma^k)\leq 
\haus^2 (\Sigma)+1/32$ for any $k>K$.
\end{itemize}
We want to prove that $\eps$ satisfies the requirement of
the lemma. Indeed choose any smooth isotopy $\varphi$ which is supported
on $B_\eps (x)$ and such that $\haus^2 (\varphi (1, \Sigma^k))\leq
\haus^2 (\Sigma^k)$.  Set
$$
K\;=\;\sup_t \haus^2 (\varphi (t, \Sigma^k)\cap B_{\sigma_k} (x))
$$
and choose $t$ such that $\mu (1-t)^2 K\leq 1/32$. 

\begin{figure}[htbp]
\begin{center}
    \input{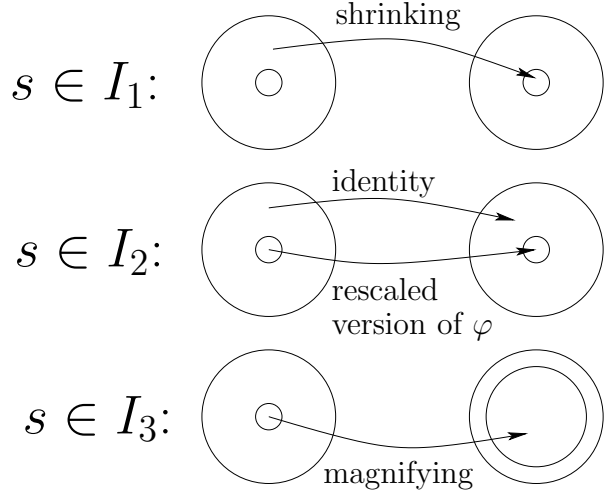}
    \caption{The isotopy $\Psi$.}
    \label{f:5}
\end{center}
\end{figure}

Define isotopies $\psi^-$ and $\tilde{\varphi}$ by
$$
\psi^-(s, \cdot)\;=\;\psi(1-s, \cdot) \qquad \text{and}\qquad 
\tilde{\varphi} \;=\; 
I_{(1-t)} \circ \varphi \circ I_{(1-t)^{-1}}\, ,
$$
and note that $\psi^-$ is the ``backward'' of $\psi$ and hence instead of 
``squeezing'', it magnifies, whereas $\tilde{\varphi}$ is
the ``$(1-t)$--shrunk'' version of $\varphi$.
We now define a (piecewise) smooth isotopy
$\Psi:I_1\cup I_2\cup I_3\times M\to M$ by gluing $\psi$, $\tilde \phi$, 
and $\psi^-$ together.  Namely, set 
\begin{itemize}
\item $\Psi(s, \cdot)=\psi(s, \cdot)$ for $s\in I_1=[0, 1-t]$;
\item $\Psi(s, \cdot)= \tilde{\varphi} (s- (1-t))$ for 
$s\in I_2=[1-t, 2-t]$;
\item $\Psi (s, \cdot)=\psi^- ((s-(2-t))+t,\tilde{\varphi} (1,\cdot))$ 
for $s\in I_3=[2-t, 3-2t]$.
 \end{itemize}
Loosely speaking, the isotopy $\Psi$ performs the following
operations (see Fig.~\ref{f:5}):
\begin{itemize}
\item[] First it ``shrinks'' the ball $B_{s_k} (x)$ to the 
ball $B_{(1-t)s_k} (x)$;
\item[] Then it applies the ``$(1-t)$--shrunk'' version of $\varphi$ 
to the ball $B_{(1-t)s_k} (x)$;
\item[] Finally it magnifies back $B_{(1-t)s_k} (x)$ 
to the ball $B_{s_k} (x)$.
\end{itemize}

By changing parameter we can assume that $\Psi$ is smooth and that 
$I_i=[(i-1)/3, i/3]$.
Note that (\ref{core}) implies that
\begin{equation}\label{ancora1}
\haus^2 (\Psi (s, \Sigma^k))\;\leq\; \haus^2 (\Sigma^k)+1/32 \qquad
\mbox{for $s\in [0,1/3]$.}
\end{equation}
Since $\mu (1-t)^2 K\leq 1/32$ we have 
\begin{equation}\label{ancora2}
\haus^2 (\Psi (s, \Sigma^k))\;\leq\; \haus^2 (\Sigma^k)+\mu
(1-t)^2K\;\leq\; \frac{1}{32}\qquad \mbox{for $s\in [1/3,2/3]$.}
\end{equation}
Again by (\ref{core}), and since $\haus^2 (\varphi (1,
\Sigma^k))\leq \haus^2 (\Sigma^k)$, we have 
\begin{eqnarray}
\haus^2 (\Psi (s, \Sigma^k))&\leq& \haus^2 (\Sigma^k) +1/32 + \mu
\haus^2 (\varphi (1, \Sigma^k)\cap B_{\sigma_k} (x))\nonumber\\
&\leq& \haus^2 (\Sigma^k) +1/32 + \mu C\eps^2\nonumber\\
&\leq& \haus^2 (\Sigma^k) +1/16 \qquad\mbox{for $s\in [2/3,1]$.}\qquad
\label{ancora3}
\end{eqnarray}
Thus for every $s$
\begin{equation}\label{lafine?}
\haus^2 (\Psi (s, \Sigma^k))\;\leq\; \haus^2 (\Sigma^k)+1/16\;\leq\;
\haus^2 (\Sigma)+3/32.
\end{equation}
Note also that $\Psi (1, \Sigma^k)= \varphi (1, \Sigma^k)$. Finally,
recall that $\Sigma^k$ was obtained via an isotopy $\varphi^k$ such
that
$$
\haus^2 (\varphi^k (t, \Sigma))\;\leq\; \haus^2 (\Sigma)+1/8.
$$
Gluing together $\varphi_k$ and $\Psi$ we easily obtain a $\Phi$
which satisfies both (\ref{itcoincides}) and (\ref{itisOK}). Clearly
the $\eps$ found in this proof satisfies the last requirement of the 
statement of the lemma.
\qed

\subsection{Proof of Theorem~\ref{goal}}\label{ss:goal} We will apply
Proposition~\ref{reg}. From Proposition~\ref{repl} we know that in
every annulus $\an \in \An_{r(x)} (x)$ 
there is a replacement $V^*$ for $V$. We still need to show that 
$V$ satisfies (a), (b), and (c) in Definition~\ref{goodrepl}.
Consider the family of surfaces $\Sigma^{j,k}$ of Lemma~\ref{stablemin}. 
By a diagonal argument we can extract a subsequence
$\Sigma^{j,k(j)}$ converging to $V^*$. 
Note the following consequence of the way we constructed 
$\{\Sigma^{j,k(j)}\}^j$. If $U$ is open and
\begin{itemize}
\item[-] either $U\cup\an$ is contained in some annulus $\An_{r(x)} (x)$
\item[-] or $U\cap \an =\emptyset$ and $U$ is contained in some
  annulus of $\An_{r(y)} (y)$ with $y\neq x$,
\end{itemize}
then $\Sigma^{j,k(j)}$ is a.m.\ in $U$. 
Thus $\{\Sigma^{j, k(j)}\}$ is still a.m.\ in 
\begin{itemize}
\item[-] every annulus of $\An_{r(x)} (x)$;
\item[-] every annulus of $\An_{\rho (y)} (y)$ for $y\neq x$, 
provided $\rho (y)$ is sufficiently small.
\end{itemize}
This shows that (b) in Definition~\ref{goodrepl} holds for $V$. Similarly, 
we can show that also condition (c) of that Definition holds.
Hence Proposition~\ref{reg} applies and we conclude that $V$ is a
smooth surface.
\qed

\appendix

\section{Proof of Proposition \ref{p:morse}}\label{a:iso}
Let $M^3$ be a closed Riemannian $3$--manifold with a Morse function
$f:M\to [0,1]$. Denote by $\Sigma_t$ the level set 
$f^{-1} (\{t\})$ and let $\Lambda$ be the saturated set of
families 
\begin{eqnarray*}
\lefteqn{\Big\{\{\Gamma_t\} \Big|\,\, \mbox{$\Gamma_t = \psi (t, \Sigma_t)$
for some $\psi\in C^\infty([0,1] \times M, M)$}}&&\\
&&\mbox{with $\psi_t \in \diff$ for every $t$}\Big\}\qquad\qquad\qquad\, .
\end{eqnarray*}

To prove Proposition~\ref{p:morse} we need to show that $m_0
(\Lambda)>0$. To do that set 
$U_t = f^{-1} ([0,t[)$ and $V_t = \psi (t, U_t)$. 
Clearly $\Gamma_t=\partial V_t$ and if we let $\Vol$ denote the volume on 
$M$, then $\Vol (U_t)$ is a
continuous function of $t$. Since $V_0$ is a finite set of points
and $V_1=M$, then there exists an $s$ such that $\Vol (V_s)= \Vol (M)/2$.
By the isoperimetric inequality there exists a constant $c(M)$ such that 
$$
\frac{\Vol (M)}{2}\;=\; 
\Vol (V_s)\;\leq\; c(M) \haus^2 (\Gamma_s)^{3/2}\, .
$$
Hence,
\begin{equation}\label{e:above}
\F (\{\Gamma_t\})\;=\;
\max_{t\in [0,1]} \haus^2 (\Gamma_t)\;\geq\; \left(
\frac{\Vol (M)}{2c(M)}\right)^{\frac{2}{3}}>0\, ,
\end{equation}
and the proposition follows.

\section{Two lemmas about varifolds}\label{App2}
\begin{lemma}\label{PittsPr} Let $U$ be an open
subset of a $3$--manifold $M$ and $W$ a stationary $2$--varifold in
$\mathcal{V} (U)$. If $K\subset\subset U$
is a smooth strictly convex set and 
$x\in (\supp (\|W\|))\cap \partial K$, then
$$
(B_r(x)\setminus \overline{K})\cap \supp (\|W\|)\;\neq\;
\emptyset\qquad \mbox{for every $r>0$.}
$$
\end{lemma}

\begin{proof} For simplicity assume that $M=\rn{3}$.
The proof can be easily adapted to the general
  case. Let us argue by contradiction; so assume that there are $x\in \supp
  (\|W\|)$ and 
 $B_r (x)$ such that $(B_r(x)\setminus \overline{K})\cap \supp (\|W\|)=
\emptyset$. Given a vector field $\chi\in C^\infty_c
  (U,\rn{3})$ and a $2$--plane $\pi$ we set
$$
\trac (D\chi (x), \pi)\;=\; D_{v_1} \chi (x)\cdot v_1 + 
D_{v_2} \chi (x)\cdot v_2
$$
where $\{v_1, v_2\}$ is an orthonormal base for $\pi$. 
Recall that the first variation of $W$ is given by
$$
\delta W (\chi)\;=\;\int_{G (U)} \trac (D\chi (x), \pi)\, dW (x, \pi)\, .
$$
Take an increasing function $\eta\in C^\infty ([0,1])$ which 
vanishes on $[3/4,1]$ and is identically 1 on $[0,1/4]$. Denote by
$\varphi$ the function given by $\varphi (x)=\eta (|y-x|/r)$ for
$y\in B_r (x)$. Take the interior unit normal $\nu$ to $\partial K$ in
$x$, and let $z_t$ be the point $x+t\nu$.  If we define vector fields 
$\psi_t$ and $\chi_t$ by 
$$
\psi_t (y)\;=\;-\frac{y-z_t}{|y-z_t|}\qquad\mbox{and}\qquad 
\chi_t= \varphi \psi_t\, ,
$$
then $\chi_t$ is supported in
$B_r (x)$ and $D \chi_t= \varphi D\psi_t +
\nabla\varphi \otimes \psi_t$. Moreover, 
by the strict convexity of the subset $K$,
$$
\nabla \varphi(y) \cdot \nu>0 \qquad \mbox{if $y\in \ov{K}\cap B_r (x)$ and
  $\nabla \varphi (y)\neq 0$.}
$$
Note that $\psi_t$ converges to $\nu$ uniformly in $B_r
(x)$, as $t\uparrow \infty$. Thus, 
$\psi_T (y)\cdot \nabla \varphi (y)\geq 0$ for 
every $y\in \ov{K}\cap B_r (x)$, provided $T$ is sufficiently large. 
This yields that 
\begin{equation}\label{positivo}
\trac
(\nabla\varphi (y) \otimes \psi_T (y), \pi)\geq 0 \qquad \mbox{for all
  } (y, \pi)\in G(B_r (x)\cap \ov {K})\, .
\end{equation}
Note that $\trac (D\psi_t (y), \pi)>0$ 
for all $(y, \pi)\in G(B_r (x))$ and all $t>0$. Thus 
\begin{eqnarray*}
\delta W (\chi_T) &=& \int_{G (B_r (x)\cap \overline{K})} 
\trac (D\chi_T (y), \pi)\, dW(y,\pi)\\
&\stackrel{(\ref{positivo})}{\geq}& 
\int_{G (B_r (x)\cap \overline{K})} \trac (\varphi (y)
D\psi_T (y), \pi)\, dW(y, \pi)\\
&\geq& \int_{G (B_{r/4} (x)\cap \overline{K})}\trac 
(D\psi_T (y), \pi)\, dW(y, \pi)\;>\; 0.\\
\end{eqnarray*}
This contradicts that $W$ is stationary and completes the proof.
\end{proof}

\begin{lemma}\label{l:tec}
Let $x\in M$ and $V$ be a stationary integer rectifiable varifold in
$M$. Assume $T$ is the subset of the support of $\|V\|$ given by
$$
T\;=\;\{\mbox{$T(y,V)$ consists of a plane 
transversal to $\partial B_{d(x,y)} (x)$}\}\,.
$$
If $\rho< \Inj (M)$, then $T$ is dense 
in $(\supp (\|V\|))\cap B_\rho (x)$.
\end{lemma}
\begin{proof} Since $V$ is integer rectifiable, then $V$ is supported on a 
rectifiable $2$--dimensional set $R$ and there exists a 
Borel function $h:R\to \N$ such that 
$V=h R$. Assume the lemma is false; 
then there exists $y\in B_\rho (x)\cap \supp (\|V\|)$ and $t>0$
such that
\begin{itemize}
\item the tangent plane to $R$ in $z$ is tangent 
to $\partial B_{d(z,x)} (x)$, for any $z\in B_t (y)$.
\end{itemize}
We choose $t$ so that $B_t (y)\subset B_\rho (x)$. Take polar coordinates 
$(r, \theta, \varphi)$ 
in $B_\rho (x)$ and let $f$ be a smooth nonnegative function 
in $C^\infty_c (B_t (y))$ with $f=1$ on $B_{t/2} (y)$.
Denote by $\chi$ the vector field 
$\chi (\theta, \varphi, r)
= f(\theta, \varphi, r) {\textstyle \frac{\partial}{\partial r}}$.
We use the notation of the proof of Lemma \ref{PittsPr}. For
every $z\in R\cap B_t (x)$, the plane $\pi$ tangent to $R$ in $z$ is
also tangent to the sphere $\partial B_{d (z,x)} (x)$.
Hence, an easy computation yields that $\trac (\chi, \pi) (z)=
2\psi (z)/d(z,x)$. This gives
$$
[\delta V] (\chi) = \int_{R\cap B_t (y)} \frac{2h(z)\psi
  (z)}{\d(z,x)} 
d\haus^2 (z) \;>\; C \|V\| (B_{t/2} (y))\, ,
$$
for some positive constant $C$. Since $y\in \supp (\|V\|)$, 
we have \
$$
\|V\| (B_{t/2} (y))\;>\;0\, .
$$ 
This contradicts that $V$ is stationary.
\end{proof}

\section{An example}\label{App1}
\begin{figure}[htbp]
\begin{center}
    \input{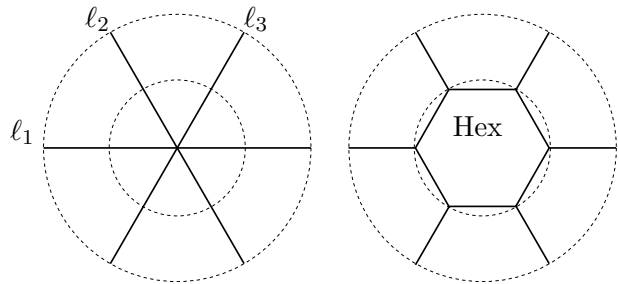}
    \caption{The varifolds $V_1$ and $V_2$.}
    \label{f:3}
\end{center}
\end{figure}

Let $V_1\in \mathcal{V}^1 (\D_2)$ be the $1$--dimensional varifold given
by three straight lines $\ell_1, \ell_2, \ell_3$ 
which meet in the origin at angles of $60$ degrees and let 
$V_2$ be the $1$--dimensional varifold given by (see Fig.~\ref{f:3}): 
\begin{itemize}
\item $V_2= V_1$ in $\D_2\setminus \D_1$;
\item In $\D_1$, $V_2$ is given by the regular hexagon ${\rm Hex}$
with sides of Length $1$ and vertices lying on the $l_i$'s. 
\end{itemize}
Note that both $V_1$ and $V_2$ are stationary in $\D_2$, they have the same
mass, and they coincide in $\D_2\setminus \D_1$.

\end{document}